\documentclass[11pt]{amsart}
\usepackage{amsxtra}
\usepackage{amssymb,stmaryrd}
\theoremstyle{plain}
\newcommand{\cleqn}{\setcounter{equation}{0}}
\newcommand{\clth}{\setcounter{theorem}{0}}
\newcommand {\sectionnew}[1]{\section{#1}\cleqn\clth}
\newtheorem{theorem}{Theorem}[section]
\newtheorem{lemma}[theorem]{Lemma}
\newtheorem{definition-theorem}[theorem]{Definition-Theorem}
\newtheorem{proposition}[theorem]{Proposition}
\newtheorem{corollary}[theorem]{Corollary}
\newtheorem{definition}[theorem]{Definition}
\newtheorem{example}[theorem]{Example}
\newtheorem{remark}[theorem]{Remark}
\newtheorem{conjecture}[theorem]{Conjecture}
\newtheorem{notation}[theorem]{Notation}
\newcommand \bth[1] { \begin{theorem}\label{t#1} }
\newcommand \ble[1] { \begin{lemma}\label{l#1} }

\newcommand \bpr[1] { \begin{proposition}\label{p#1} }
\newcommand \bco[1] { \begin{corollary}\label{c#1} }
\newcommand \bde[1] { \begin{definition}\label{d#1}\rm }
\newcommand \bex[1] { \begin{example}\label{e#1}\rm }
\newcommand \bre[1] { \begin{remark}\label{r#1}\rm }
\newcommand \bcj[1] { \begin{conjecture}\label{j#1}\rm }

\newcommand \bnota[1] { \begin{notation}\label{n#1}\rm }
\renewcommand {\eth} { \end{theorem} }
\newcommand {\ele} { \end{lemma} }

\newcommand {\epr} { \end{proposition} }
\newcommand {\eco} { \end{corollary} }
\newcommand {\ede} { \end{definition} }
\newcommand {\eex} { \end{example} }
\newcommand {\ere} { \end{remark} }
\newcommand {\ecj} { \end{conjecture} }

\newcommand {\enota} { \end{notation} }
\newcommand \thref[1]{Theorem \ref{t#1}}
\newcommand \leref[1]{Lemma \ref{l#1}}
\newcommand \prref[1]{Proposition \ref{p#1}}
\newcommand \coref[1]{Corollary \ref{c#1}}

\newcommand \deref[1]{Definition \ref{d#1}}

\newcommand \reref[1]{Remark \ref{r#1}}

\def \Rset {{\mathbb R}}         %mathsets
\def \Cset {{\mathbb C}}
\def \KK {{\mathbb K}}

\def \Zset {{\mathbb Z}}
\def \Nset {{\mathbb N}}
\def \Qset {{\mathbb Q}}

\newcommand\HH{{\mathcal{H}}}     %mathcal
\def \CC {{\mathcal{C}}}
\def \DD {{\mathcal{CD}}}
\def \CG {{\mathcal{CG}}}
\def \LP {{\mathcal{LP}}}
\def \RP {{\mathcal{RP}}}
\def \II {{\mathcal{I}}}

\def \UU {{\mathcal{U}}}
\def \RR {{\mathcal{R}}}
\def \SS {{\mathcal{S}}}

\def \TT {{\mathcal{T}}} 

\def \De {\Delta}   % Greek letters
\def \de {\delta}
\def \al {\alpha}
\def \be {\beta}
\def \vpi {\varpi}
\def \la {\lambda}

\def \La {\Lambda}

\def \om {\omega}

\def \ga {\gamma}
\def \de {\delta}

\def \Sig {\Sigma}
\def \sig {\sigma}

\def \vp {\varphi}

\def \sig{\sigma}

\def \pp {p}
\def \mt  {\mapsto}
\def \hra {\hookrightarrow}

\def \st  {\ast}                 %duals
\def \ci  {\circ}

\def \rcor {\rangle}
\def \lcor {\langle}
\newcommand \ola[1] {#1^\text{\tiny \hspace{-10pt} $\leftarrow$}}
\def \ol {\overline}
\def \wt {\widetilde}

\def \id { {\mathrm{id}} }

\def \Lie { {\mathrm{Lie \,}} }

\def \g  {\mathfrak{g}}   % Lie algebra letters

\def \sl {\mathfrak{sl}}

\def \n  {\mathfrak{n}}

\def \b  {\mathfrak{b}}

\def \sl {\mathfrak{sl}}
\DeclareMathOperator \Ker { {\mathrm{Ker}} }

\DeclareMathOperator \GKdim {{\mathrm{GK \, dim}}}

\DeclareMathOperator \lt  { {\mathrm{lc}} }
\DeclareMathOperator \rev {{\mathrm{rev}} }

\DeclareMathOperator \Fract { {\mathrm{Fract}} }

\renewcommand \max { {\mathrm{max}} }
\newcommand \Spec { {\mathrm{Spec}} }

\begin{document}
%%%%%%%%%%%%%%%%%%%%%%%%%%%%%%%%%%%%%%%%%%%%%%%%%%%%%%%%%%%%%%%%%%%%%%%%%%%
%%%%%%%%%%%%%%%%%%%%%%    Title    %%%%%%%%%%%%%%%%%%%%%%%%%%%%%%%%%%%%%%%%
\title[Quantum Schubert cell algebras and Richardson varieties]
{Prime factors of quantum Schubert cell algebras
and clusters for quantum Richardson varieties}
\author[T. H. Lenagan]{T. H. Lenagan}
\address{
Maxwell Institute for Mathematical Sciences \\
School of Mathematics   \\
University of Edinburgh \\
James Clerk Maxwell Building \\
King's Buildings \\
Peter Guthrie Tait Road \\
Edinburgh EH9 3FD \\
Scotland \\
United Kingdom
}
\email{T.Lenagan@ed.ac.uk}
\author[M. T. Yakimov]{M. T. Yakimov}
\thanks{The research of T.H.L. was partially supported by EPSRC grant EP/K035827/1 and 
a Leverhulme Trust Emeritus Fellowship. The research of M.T.Y. was partially supported by NSF grant 
DMS-1303038 and Louisiana Board of Regents grant Pfund-403}
\address{
Department of Mathematics \\
Louisiana State University \\
Baton Rouge, LA 70803 \\
U.S.A.
}
\email{yakimov@math.lsu.edu}
\date{}
\begin{abstract} The understanding of the topology of the spectra of quantum Schubert 
cell algebras hinges on the description of their prime factors by ideals invariant under the 
maximal torus of the ambient Kac--Moody group. We give an explicit description of these prime quotients
by expressing their Cauchon generators 
in terms of sequences of normal elements in chains of subalgebras. Based on this, we construct large families of 
quantum clusters for all of these algebras and the quantum Richardson varieties associated to arbitrary 
symmetrizable Kac--Moody algebras and all pairs of Weyl group elements. Along the way we develop a quantum version 
of the Fomin--Zelevinsky twist map for all quantum Richardson varieties. Furthermore, we establish an 
explicit relationship between the Goodearl--Letzter and Cauchon approaches to the descriptions of the 
spectra of symmetric CGL extensions. 
\end{abstract}
\maketitle
%%%%%%%%%%%%%%%%%%%%   Introduction   %%%%%%%%%%%%%%%%%%%%%%%%%%%%%%%%%%%%%%%%
\sectionnew{Introduction}
\label{intro}
\subsection{Background}
\label{1.1}
The quantum Schubert cell algebras play an important role in representation theory (the Kashiwara--Lusztig theory 
of crystal/canonical bases \cite{K0,L0}), ring theory \cite{J,MC,Y-plms}, Hopf algebras (coideal subalgebras \cite{HS}) and 
cluster algebras \cite{GLS,GY2}. Let $\g$ be a symmetrizable Kac--Moody algebra and $w$ an element of its 
Weyl group. The corresponding quantum Schubert cell algebras $\UU^\pm[w]$ are deformations of the 
universal enveloping algebras $\UU(\n_\pm \cap w(\n_\mp))$ where $\n_\pm$ are the nilradicals of the 
standard opposite Borel subalgebras of $\g$. They were defined 
by De Concini--Kac--Procesi \cite{DKP},  Lusztig \cite{L} and by Beck in the affine case \cite{B}.

In this paper we construct explicit models for the prime quotients of the quantum Schubert cell algebras
by ideals invariant under the maximal torus of the ambient Kac--Moody group.
These quotients play a key role in two problems that have attracted a lot of attention. 
One is of algebraic nature and is about the description of the topology
of the spectra of quantum Schubert cell algebras. The other is of combinatorial nature -- the construction of 
cluster algebra structures on quantum and classical Richardson varieties.   

From the point of view of ring theory, the algebras $\UU^-[w]$ are large families of 
deformations of universal enveloping algebras of nilpotent Lie algebras. 
It is a long-standing problem to carry out an analog of the orbit method \cite{Di}
for these types of algebras. 
The canonical maximal torus $\HH$ of the related Kac--Moody group $G$ 
acts on $\UU^-[w]$ by algebra automorphisms. The $\HH$-invariant prime ideals of $\UU^-[w]$ were classified in 
\cite{MC,Y-plms}, where it was shown that they are parametrized by $W^{\leq w}:= \{ u \in W \mid u \leq w\}$ -- 
the corresponding ideals will be denoted by $I_w(u)$. By a general result of Goodearl and Letzter \cite{GL}, 
$\Spec \, \UU^-[w]$ is partitioned into 
\[
\Spec \, \UU^-[w] = \bigsqcup_{u \in W^{\leq w}} \Spec_u \, \UU^-[w]
\]
where each stratum $\Spec_u \, \UU^-[w]$ is homeomorphic to a torus and the ideals 
in it are obtained by extension and contraction from the center of $\Fract(\UU^-[w]/I_w(u))$. 
The main ring-theoretic problem for $\UU^-[w]$ is to describe the topology of their spectra, 
ideally by identifying it 
with the topological space of the symplectic foliation of the standard Poisson structure on the 
full flag variety of $G$ (restricted to the Schubert cell corresponding to $w$). Understanding this topology 
amounts to solving the containment problem for the prime ideals of $\UU^-[w]$ for which one needs an 
explicit model for the $\HH$-prime quotients $\UU^-[w]/I_w(u)$. 

Recall that a quantum torus is an algebra of the form 
\[
\TT := \frac{\KK \lcor Y_1^{\pm 1}, \ldots, Y_n^{\pm 1} \rcor}{(Y_k Y_j - q_{kj} Y_j Y_k, k>j )} 
\]
for some $q_{kj} \in \KK^*$. From each reduced expression 
\begin{equation}
\label{red-expr0}
w = s_{i_1} \ldots s_{i_N}
\end{equation}
one constructs Lusztig's root vectors $F_{\be_1}, \ldots, F_{\be_N}$ which form a generating set
of $\UU^-[w]$; here $\be_1, \ldots, \be_N$ are the roots of $n_+ \cap w(\n_-)$. 
The algebra $\UU^-[w]$ has two presentations 
as an iterated skew polynomial extension
\begin{align}
\label{pres1}
\UU^-[w] &= \KK [F_{\be_1}] [F_{\be_2}; \sig_2, \de_2] \ldots [F_{\be_N}; \sig_N, \de_N]
\\
         &= \KK [F_{\be_N}] [F_{\be_{N-1}}; \sig^*_{N-1}, \de^*_{N-1}] \ldots [F_{\be_1}; \sig^*_1, \de^*_1].
\label{pres2}
\end{align}
The Cauchon method of deleted derivations (applied to the first presentation) 
constructs in an iterative fashion a set of elements 
$Y_k \in \Fract (\UU^-[w]/I_w(u))$, indexed by a subset $D(u) \subseteq [1,N]$ 
called the Cauchon diagram of $I_w(u)$, such that the elements $\{Y_k^{\pm 1} \mid k \in D(u) \}$ 
generate a copy of a quantum torus 
\[
\TT(Y_k , k \in D(u)) \hra \Fract (\UU^-[w]/I_w(u))
\]
satisfying
\[
\UU^-[w]/I_w(u) \hra \TT(Y_k , k \in D(u)).
\]
This is the key construction that is currently used to understand the factors $\UU^-[w]/I_w(u)$ and 
the topology of $\Spec \, \UU^-[w]$. 
There are two main difficulties with it. Firstly,
the elements $Y_k$ are the result of an involved iterative 
construction and are not explicit in any way.
Secondly, the elements $Y_k$ do not lie 
in the algebra $\UU^-[w]/I_w(u)$ in general but rather in its division ring of fractions. This leads to difficulties for
the Goodearl--Letzter approach to $\Spec \, \UU^-[w]$ because one needs to contract ideals 
of $\TT(Y_k , k \in D(u))$ to $\UU^-[w]/I_w(u)$.  

The Cauchon method can be also applied to the reverse presentation \eqref{pres2} giving rise to 
another quantum torus in which $\UU^-[w]/I_w(u)$ is embedded 
\[
\UU^-[w]/I_w(u) \hra \TT(Y_{k, \rev}, k \in D_{\rev} (u)) \hra \Fract (\UU^-[w]/I_w(u)).
\]
The are analogous difficulties in this situation.
%%%%%%%%%%%%%%%%%%
\subsection{Results on the description of Cauchon generators}
\label{1.2}
We resolve the above problems and show that the Cauchon tori $\TT(Y_k , k \in D(u))$ 
and $\TT(Y_{k, \rev} , k \in D_{\rev}(u))$
have explicit generating sets that lie in $\UU^-[w]/I_w(u)$.
Our approach is based on the following general idea. Denote the canonical projection 
\[
p \colon \UU^-[w] \to \UU^-[w]/I_w(u).
\]
Consider chains $\CC$ of subalgebras
\[
A_1 \subset A_2 \subset \ldots \subset A_N = \UU^-[w]
\]
such that
\begin{equation}
\label{GKprop}
\GKdim A_k - \GKdim A_{k-1} =1.
\end{equation}
Here and below $\GKdim$ denotes the Gelfand--Kirillov dimension of an 
algebra. Given such a chain, we can project it
\[
\pp(A_1) \subset \pp(A_2) \subset \ldots \subset \pp(A_N) = 
\UU^-[w]/I_w(u)
\]
and associate to it a subset of the form
\begin{multline*}
\Sig_\CC:=
\{\mbox{a nonzero normal element $y_k \in \pp(A_k)$} 
\\
\mbox{for those $k$ such that 
$\GKdim \pp(A_k) - \GKdim \pp(A_{k-1}) =1$} \}.
\end{multline*}
We will call the set of $k$'s, the {\em{jump set}} of the chain $\CC$. 
The main idea is to construct chains with the property that the projected subalgebras have sufficiently 
many nontrivial normal elements, and to show that sets $\Sig_\CC$ of the above form 
generate quantum tori and $\UU^-[w]/I_w(u)$ embeds in them.

For each reduced expression \eqref{red-expr0}, there are two special chains of subalgebras of $\UU^-[w]$ 
obtained by adjoining the generators \eqref{pres1} of $\UU^-[w]$ in the direct order and the 
reverse order \eqref{pres2}. Denote them by $\CC$ and $\CC_{\rev}$. Our first result is as follows, 
see Theorems \ref{tmain} and \ref{tmain2} for details.
\medskip
\\
\noindent
{\bf{Theorem A.}} {\em{(i) The jump set of the chain $\CC$ is the complement to 
the index set of the unique right positive subexpression of 
\eqref{red-expr0} with total product $u$ in the sense of Deodhar, Marsh and Rietsch \cite{De,MR}. A sequence of normal elements 
$\Sig_\CC$ is provided by the sequence of quantum minors 
\[
\pp(\De_{u_{\leq k} \vpi_k, w_{\leq k} \vpi_k}) 
\]
for the integers $k$ in this index set. The set $\Sig_\CC$ is a set of independent generators of 
the Cauchon quantum torus $\TT(Y_{k, \rev} , k \in D_{\rev}(u))$ and these generators belong to $\UU^-[w]/I_w(u)$. 
Furthermore, the Cauchon generators $Y_{k, \rev}$ are Laurent monomials in the projected quantum minors 
whose exponents form a triangular matrix.

(ii) The jump set of the chain $\CC_{\rev}$ is the complement to the index set of the unique left positive subexpression of 
\eqref{red-expr0} with total product $u$. A sequence of normal elements $\Sig_{\CC_{\rev}}$ can be constructed
using a similar sequence of quantum minors; $\Sig_{\CC_{\rev}} \subset \UU^-[w]/I_w(u)$ is a set of independent generators of 
the Cauchon quantum torus $\TT(Y_k, k \in D(u))$. The Cauchon generators $Y_k$ are Laurent monomials in 
the elements of $\Sig_{\CC_{\rev}}$ with exponents forming a triangular matrix. 
}}
\medskip
\\
In Theorems \ref{tUw-contr} and \ref{tUw-contr2} we also describe explicitly all projected algebras $\pp(A_k)$ and 
$\pp(A_{k, \rev})$, in other words all contractions
\[
I_w(u) \cap A_k \quad \mbox{and} \quad I_w(u) \cap A_{k, \rev}.
\]    
%%%%%%%%%%%%%%%%%%
\subsection{Results on quantum twists maps for Richardson varieties in symmetrizable Kac--Moody 
groups}
\label{1.3}
The full flag variety of a symmetrizable Kac--Moody group $G$ has the   
Schubert cell decompositions
\[
G/B_+ = \bigsqcup_{w \in W} B_+ \cdot w B_+ = \bigsqcup_{u \in W} B_- \cdot u B_+
\]
where $B_\pm$ is a pair of opposite Borel subgroups. 
The open Richardson varieties are defined by 
\[
R_{u,w} = B_+ \cdot w B_+ \cap B_- \cdot u B_+,
\]
and $R_{u, w} \neq \varnothing$ if and only if $u \leq w$, \cite[Corollary 1.2]{De}. Since
$\Lie(B_+) + \Lie (B-) = \Lie (G)$, 
they are smooth, irreducible varieties 
and $\dim R_{u,w}= \ell(w) - \ell(u)$. We have 
\[ 
G/B_+ = \bigsqcup_{u \leq w \in W} R_{u,w}.
\]
The quantized coordinate ring of $R_{u,w}$ can be expressed as 
\[
R_q[R_{u,w}]:= (\UU^-[w]/I_w(u)) [E_{u,w}^{-1}]
\]
where $E_{u,w}$ is a multiplicative subset of $\UU^-[w]/I_w(u)$ 
consisting of normal elements, see \S \ref{2.5} for details.

Twist maps, defined by Fomin and Zelevinsky \cite{FZ}, are certain isomorphisms 
between double Bruhat cells. They play a major role in the study of the totally nonnegative 
part of $G$, canonical bases and cluster structures for double Bruhat cells. More recently, such were considered 
for open Richardson varieties in Grassmannians \cite{MS,MSp}, and were used to study cluster expansions 
and to prove local acyclicity for the related cluster algebras. We construct twist maps 
for all Richardson varieties and, in addition, do this in the quantum situation.
\medskip
\\
\noindent
{\bf{Theorem B.}} {\em{There is an algebra antiisomorphism
\[
\Theta_w \colon \UU^-[w^{-1}] \to \UU^-[w]
\]
given by \eqref{twist}. It satisfies $\Theta_w(I_{u^{-1}}(w^{-1})) = I_u(w)$
for all $u \in W$, $u \leq w$ and induces an antiisomorphism
\[
\Theta_w \colon R_q[R_{u^{-1}, w^{-1}}] \to R_q[R_{u,w}].
\]
}}
The twist map $\Theta_w$ interchanges the statements in parts (i) and (ii) of Theorem A.
%One easily construct a classical counterpart of $\Theta_w$ from the semiclassical limit.
%%%%%%%%%%%%%%%%%%
\subsection{Results on quantum clusters for Richardson varieties in symmetrizable Kac--Moody 
groups}
Recently, for each symmetric Kac--Moody algebra $\g$, Leclerc \cite{Le} defined a cluster algebra 
structure inside the coordinate ring of each Richardson variety $R_{u,w}$ such that the two algebras 
have the same dimension. We apply Theorem A to obtain large families of toric frames 
for the algebras $\UU^-[w]/I_w(u)$ and $R_q[R_{u,w}]$, with the ultimate goal of controlling 
the size of Leclerc's cluster algebra from below.
Similarly to \cite{GY2}, consider the following subset of the symmetric group $S_N$:
\[ 
\Xi_N := \{ \pi \in S_N \mid \pi([1,k]) \; \; \mbox{is an interval for all} \; \; k \in [2,N] \}.
\]
The chain of subalgebras obtained by adjoining the Lusztig generators 
of $\UU^-[w]$ in the order $F_{\be_{\pi(1)}}, \ldots, F_{\be_{\pi(N)}}$ 
has the property \eqref{GKprop}. Denote this chain by $\CC_\pi$. We recover 
the chains $\CC$ and $\CC_{\rev}$ for the identity and the longest element of $S_N$.
Those are elements of the set $\Xi_N$ which is very large.
\medskip
\\
\noindent
{\bf{Theorem C.}} {\em{For all elements $\pi \in \Xi_N$, one can construct sets 
$\Sig_{\CC_\pi}$ consisting of projected quantum minors. Each of these sets produces 
toric frames for $\UU^-[w]/I_w(u)$ and $R_q[R_{u,w}]$ in the sense of Berenstein and Zelevinsky, 
\cite{BZ}.}} 
\medskip
\\
\noindent
For the detailed formulation of this result we refer to \thref{toric-frames}. 
%%%%%%%%%%%%%%%%%%
\subsection{Unifying the Cauchon and Goodearl--Letzter approaches to the torus invariant prime ideals 
of CGL extensions}
\label{1.5}
The quantum Schubert cell algebras $\UU^-[w]$ are members of the large axiomatic class of
Cauchon--Goodearl--Letzter (CGL) extensions. These are iterated skew polynomial extensions with an action 
of a torus $\HH$ that satisfy certain natural properties resembling the definition 
of (universal enveloping algebras of) nilpotent Lie algebras via derived series, see \deref{CGL}. 
There are two approaches to describing the $\HH$-invariant prime ideals of such algebras $R$. 
The Goodearl--Letzter \cite{GL} one describes these ideals via recursive contractions with the subalgebras of 
$R$ and shows that at each step at most 2 ideals lead to the same contraction. The Cauchon approach first checks 
if a generator $x$ of $R$ belongs to an $\HH$-prime ideal $I$ and then maps the ideal to the leading coefficients of its elements 
written as polynomials in $x$, or to another contraction ideal. No connection between the two approaches 
was previously found.

In \thref{CGLcontr} we unify the two approaches for symmetric CGL extensions -- extensions that satisfy 
the CGL axioms for the direct and reverse order of adjoining the generators of $R$. This relation interchanges the 
two approaches applied to the two opposite presentations.

The results of the paper have applications to the problems in \S \ref{1.1} that will be described 
in forthcoming publications. Firstly, we will show that the toric frames in Theorem C are related 
by mutations and use this to control the size of Leclerc's cluster algebras \cite{Le} from below.
Secondly, we use Theorem A to set up a torus equivariant map from the symplectic foliation of a Schubert cell 
to the primitive spectrum of the corresponding quantum Schubert cell algebra. This will be a conjectural candidate for 
the desired homeomorphism from \S \ref{1.1} for all quantum Schubert cell algebras $\UU^-[w]$. We believe that 
this will provide a framework in which one can attempt to settle 
the Brown--Goodearl conjecture \cite[Conjecture 3.11]{BG2} in the case of the 
algebras $\UU^-[w]$; this is a general conjecture on the topology of spectra of quantum algebras that is only 
verified in very low GK dimension. Finally, we will also construct a direct relationship between the 
spectra of the quantum Schubert cell algebras $\UU^-[w]$ and the totally nonnegative part of the corresponding 
Schubert cell; previously such was obtained for the algebras of quantum matrices \cite{GLL}.
%%%%%%%%%%%%%%%%%%%%%%%%%%%%%%%%%%%%%%%%%%%%%%%%%%%%%%%%%%%%%%%%%%%%%%%%%%%%%%
\medskip
\\
\noindent
{\bf Acknowledgements.} We are thankful to Ken Goodearl, Ryan Kinser, Allen Knutson, 
St\'ephane Launois, Bernard Leclerc and Jiang-Hua Lu for helpful discussions and for their comments on 
the first draft of the paper. We would also like to thank the anonymous referee for the very helpful suggestions.
%%%%%%%%%%%%%%%%%%%%%%%%%%%%%%%%%%%%%%%%%%%%%%%%%%%%%
\sectionnew{Quantum Schubert cells, CGL extensions, 
and torus invariant primes}
\label{qalg}
%%%%%%%%%%%%%%%%
In this section we collect some facts on quantum groups and quantum Schubert cell algebras, 
as well as facts on their prime spectra, that will be used in the paper. For more details on quantum 
groups we refer the reader to \cite{Ja,J}. 
%%%%%%%%%
\subsection{Quantum algebras}
\label{2.1}
Let $\g$ be a symmetrizable Kac--Moody algebra of rank $r$ with Weyl group $W$ and 
set of simple roots $\al_i$, $i \in [1,r]$. Let
$\lcor.,. \rcor$ be the invariant bilinear form on $\Rset \al_1 + \cdots + \Rset \al_r$
normalized by $\lcor \al_i, \al_i \rcor = 2$ for short roots 
$\al_i$. Denote by $P^+$ the set of dominant integral weights of $\g$, and by 
$P$ and $Q$ the weight and root lattices of $\g$. Let $\{\vpi_i\}$ and $\{\al_i\spcheck\}$ 
be the fundamental weights and simple coroots of $\g$. The corresponding 
simple reflections in $W$ will be denoted by $\{s_i\}$. 

Let $\KK$ be an arbitrary infinite base field and $q \in \KK^*$ be a 
non-root of unity. Denote by $\UU_q(\g)$ the quantized universal enveloping 
algebra of $\g$ over the base field $\KK$ with deformation parameter $q$. 
We will use the conventions of \cite{Ja} for the 
(Hopf) algebra structure on $\UU_q(\g)$, with the exception that the generators of $\UU_q(\g)$ 
will be denoted by $E_i, F_i, K_i^{\pm 1}$, indexed by $i \in [1,r]$ rather than by the 
set of simple roots of $\g$. 
Recall that the weight spaces of a $\UU_q(\g)$-module $V$ are defined by 
\[
V_\nu := \{ v \in V \mid K_i v = q^{ \lcor \al_i, \nu \rcor} v \}, \quad \nu \in P. 
\]
For $\la \in P^+$ denote by $V(\la)$ the unique irreducible highest 
weight $\UU_q(\g)$-module with highest weight $\la$. 
Let $v_\la$ be a highest weight vector of $V(\la)$. We will use Lusztig's 
actions of the braid group of $\g$ on $\UU_q(\g)$ and $V(\la)$, $\la \in P^+$, 
in the conventions of \cite{Ja}.  

Denote by $\UU^\pm_q(\g)$ the unital subalgebras of $\UU_q(\g)$ 
generated by $\{E_i\}$ and $\{F_i\}$, respectively.
Given a Weyl group element $w$ and a reduced expression 
\begin{equation}
\label{red-expr}
w = s_{i_1} \ldots s_{i_N}
\end{equation}
of $w$, consider the root vectors
\[
\be_1 = \al_{i_1}, \be_2 = s_{i_1} (\al_{i_2}), \ldots, 
\be_N = s_{i_1} \ldots s_{i_{N-1}}(\al_{i_N}).
\]
De Concini, Kac, and Procesi \cite{DKP}, Lusztig \cite[\S 40.2]{L} and Beck \cite{B}
defined the quantum Schubert cell algebras 
$\UU^\pm[w]$ as the unital subalgebras of $\UU_q^\pm(\g)$ with generators 
\begin{multline}
\{ E_{\be_j}:=
T_{i_1} \ldots T_{i_{j-1}} 
(E_{\al_j}) \mid j \in [1,N] \} \quad 
\mbox{and} 
\\
\{ F_{\be_j} := 
T_{i_1} \ldots T_{i_{j-1}} 
(F_{\al_j}) \mid 
j \in [1, N] \},
\label{rootv}
\end{multline}
respectively, and proved that these algebras do not depend on the 
choice of a reduced expression of $w$. 
Define the quantum $R$-matrix associated to $w$ by
\begin{equation}
\label{Rw} 
\RR^w := 
\sum_{m_1, \ldots, m_N \in \Nset}
\left( \prod_{j=1}^N 
\frac{ (q_{i_j}^{-1} - q_{i_j})^{m_j}}{
{q_{i_j}}^{m_j (m_j-1)/2} [m_j]_{q_{i_j}}! } \right)
E_{\be_N}^{m_N} \ldots E_{\be_1}^{m_1} \otimes 
F_{\be_N}^{m_N} \ldots F_{\be_1}^{m_1} 
\end{equation}
considered as an element of the completion of $\UU^+ \otimes \UU^-$ 
with respect to the descending filtration \cite[\S 4.1.1]{L}. As usual,
$q$-integers and $q$-factorials are defined by 
\[
[n]_q := \frac{q^n - q^{-n}}{q- q^{-1}}, \; 
[n]_q ! := [1]_q \ldots [n]_q, \; 
n \in \Nset.
\]

For $\la \in P^+$ and $w \in W$ set 
\[
v_{w \la} := T_{w^{-1}}^{-1} v_\la \in V(\la)_{w \la}.
\]
It is well known that $v_{w \la}$ depends only on $w \la$ and not on the 
choice of $w$ and $\la$. Since $\dim V(\la)_{w \la}=1$, 
there is a unique dual vector
\[
\xi_{w \la} \in (V(\la)^*)_{- w \la} \quad
\mbox{such that} \quad
\lcor \xi_{w \la}, v_{w \la} \rcor =1.
\] 
For a pair of Weyl group elements $(u,w)$ one defines the quantum minor 
\[
c_{u \la, w \la} \in (\UU_q(\g))^* \quad
\mbox{given by} \quad
c_{u \la, w \la}(x) = \lcor \xi_{u \la}, x v_{w \la} \rcor, \; \; 
\forall x \in \UU_q(\g).
\] 
This quantum minor
does not depend on the choice of a highest weight vector of $V(\la)$. 
Given a reduced expression \eqref{red-expr} and $k \in [1, N]$, set
$w_{\leq k} := s_{i_1} \ldots s_{i_k}$. The algebras 
$\UU^\pm[w_{\leq k}]$ coincide with the subalgebras of $\UU^\pm[w]$ 
generated by the subsets of \eqref{rootv} with $j \in [1,k]$.  
For $u \in W$ and $k \in [1,N]$, consider the quantum minors 
\begin{equation}
\label{q-min}
\De_{u \la, w_{\leq k} \la} = 
\lcor c_{u \la, w_{\leq k}  \la} \tau \otimes \id, \RR^{w_{\leq k}} \rcor
\in \UU^-[w_{\leq k}] \subset \UU^-[w]
\end{equation}
where $\tau$ denotes the unique graded algebra 
antiautomorphism of $\UU_q(\g)$ defined via
\begin{equation}
\label{tau}
\tau(E_i) = E_i,
\,
\tau(F_i) = F_i, 
\, 
\tau(K_i) = K_i^{-1}, \quad
\forall i \in [1,r],
\end{equation}
see \cite[Lemma 4.6(b)]{Ja}.
The quantum minor \eqref{q-min} is nonzero iff $u \leq w_{\leq k}$ in the Bruhat order.
Using the form of $\RR^w$ and the highest weight property of $v_\la$, 
one easily derives that 
\begin{equation}
\label{q-min2}
\De_{u \la, w_{\leq k} \la} = 
\lcor c_{u \la, w_{\leq k} \la} \tau \otimes \id, \RR^{w} \rcor.
\end{equation}

\subsection{CGL extensions}
\label{2.2}
Consider an iterated skew polynomial extension of length $N$,
\begin{equation} 
\label{itOre}
R := \KK[x_1][x_2; \sig_2, \delta_2] \cdots [x_N; \sig_N, \delta_N].
\end{equation}
For $k \in [0,N]$, denote the $k$-th algebra in the 
chain of extensions 
\[
R_k:=\KK[x_1][x_2; \sig_2, \de_2] 
\cdots [x_k; \sig_k, \de_k]. 
\]
In particular, $R_0 = \KK$ and 
$R_N = R$.

\bde{CGL} An iterated skew polynomial extension $R$ as in \eqref{itOre} 
is called a \emph{Cauchon--Goodearl--Letzter} (\emph{CGL}) \emph{extension} 
if it is 
equipped with a rational action of a $\KK$-torus $\HH$ 
by $\KK$-algebra automorphisms satisfying the following conditions:
\begin{enumerate}
\item[(i)] The elements $x_1, \ldots, x_N$ are $\HH$-eigenvectors.
\item[(ii)] For every $k \in [2,N]$, $\de_k$ is a locally nilpotent 
$\sig_k$-derivation of $R_{k-1}$. 

\item[(iii)] For every $k \in [1,N]$, there exists $h_k \in \HH$ such that 
$\sig_k = (h_k \cdot)$ and the eigenvalue of $x_k$, to be denoted by $\la_k \in \KK^*$,
is not a root of unity.
\end{enumerate}
\ede

A CGL extension $R$ possesses the following canonical chain 
of subalgebras which are CGL extensions:
\begin{equation}
\label{chain}
R_1 \subset R_2 \subset \ldots \subset R_N= R,
\end{equation}
where $R_k$ are equipped with the restriction of the $\HH$-action.

For $1 \leq j < k \leq N$ denote the eigenvalues $\la_{kj} \in \KK$ given by 
\[
\sig_k(x_j) = h_k \cdot x_j = \la_{kj} x_j.
\]
Set $\la_{kk} =1$ and $\la_{jk} = \la_{kj}^{-1}$ for $j>k$. 
For $j, k  \in [1,N]$ denote by $R_{[j,k]}$ the unital subalgebra of
$R$ generated by $\{ x_i \mid j \le i \le k \}$. In particular, 
$R_{[j,k]} = \KK$ if $j \nleq k$.

\bde{sCGL} A CGL extension $R$ of length $N$ as above is called {\em{symmetric}} if it can be presented 
as an iterated skew polynomial extension for the reverse order of
its generators,
\[
R = \KK[x_N][x_{N-1}; \sig^*_{N-1}, \delta^*_{N-1}] \cdots [x_1; \sig^*_1, \delta^*_1],
\]
in such a way that conditions {\rm(ii)--(iii)} in \deref{CGL}
are satisfied for some choice of $h_N^*, \ldots, h_{1}^* \in \HH$.
\ede
The following proposition is easy to prove and is left to the reader.
\bpr{sCGL2} A CGL extension $R$ as above is symmetric if and only if it satisfies 
the Levendorskii--Soibelman type straightening law
\[
x_k x_j - \la_{kj} x_j x_k \in R_{[j,k]}, \quad 
\forall j <k 
\] 
and there exist $h_j^* \in \HH$, $\forall j \in [1,N]$, 
such that $h_j^* \cdot x_k = \la_{kj}^{-1} x_k$
for all $k>j$. In this case, the endomorphisms 
$\sig^*_j$ and $\de^*_j$ of $R_{[j+1, N]}$ 
are given by $\sig_k := (h^*_k \cdot)$ and 
\[
\de^*_j(x_k): = x_j x_k - \la_{jk} x_k x_j = - \la_{jk} \de_k(x_j), \; \; 
\forall k \in [j+1, N].
\] 
\epr
A symmetric CGL extension $R$ possesses the following reverse chain 
of subalgebras which are CGL extensions:
\begin{equation}
\label{chain-rev}
R_{N, \rev}  \subset R_{N-1, \rev} \subset \ldots \subset R_{1, \rev} = R,
\end{equation}
where the intermediate subalgebras $R_{k, \rev}$ are given by
\[
R_{k, \rev} = \KK[x_N][x_{N-1}; \sig^*_{N-1}, \delta^*_{N-1}] \cdots [x_k; \sig^*_k, \delta^*_k]
\]
and are equipped with the restriction of the $\HH$-action. 

\subsection{Cauchon's method of deleting derivations}
\label{2.3}
Consider a CGL extension $R$ as above. The Cauchon map 
\[
\theta_{x_N} \colon R_{N-1} \to R_N[x_N^{-1}] \; \; 
\mbox{is given by} \; \; 
\theta_{x_N}(b) = \sum_{m=0}^\infty 
\frac{(1- \la_N)^{-m}}{(m)_{\la_N}!} [ \delta_N^m \sig^{-m}_N(b)] 
x_N^{-m}
\]
for $b \in R_{N-1}$. We set $(0)_q=1$, $(m)_q = (1-q^m)/(1-q)$ for $m > 0$, and 
$(m)_q! = (0)_q \ldots (m)_q$ for $m \in \Nset$ and a non-root of unity $q \in \KK^*$. 
This map is an injective $\HH$-equivariant algebra homomorphism, \cite{Ca1}. 
Denote $R'_{N-1}:=\theta_{x_N}(R_{N-1})$, and let $R'_N$ be the subalgebra of $R[x_N^{-1}]=R_N[x_N^{-1}]$
generated by $R'_{N-1}$ and $x_N$. The map $\sig_N$ extends to an 
automorphism of $R_N[x_N^{-1}]$ by setting $\sig_N = (h_N \cdot)$ 
(in the notation of \deref{CGL}) because $x_N$ is an $\HH$-eigenvector, and furthermore 
$\sig_N$ restricts to an automorphism of $R'_{N-1}$. Then, \cite{Ca1}
\begin{equation}
\label{delete-der1}
R'_N \cong R'_{N-1}[x_N; \sig_N] \; \; \mbox{and} \; \; 
R_N[x_N^{-1}] = R'_N[x_N^{-1}].
\end{equation}

For an element $a = b_n x_N^n + \cdots + b_m x_N^m \in R_N[x_N^{-1}]$ with $n \leq m$ and 
$b_m\neq 0$, denote its leading coefficient
\[
\lt_{x_N}(a) := b_m
\]
(called leading term in \cite{GeY}).
Set $\lt_{x_N}(0):=0$. For a subset $S$ of $R_N[x_N^{-1}]$ denote by $\lt_{x_N}(S)$ the set of leading 
coefficients of all elements of $S$.

Goodearl and Letzter proved \cite[Proposition 4.2]{GL} all $\HH$-prime ideals of 
a CGL extension are completely prime. 
Cauchon's method of deleting derivations associates to each 
$\HH$-prime ideal $I$ of a CGL extension $R$ of length $N$ as above, 
a subset $\DD(I) \subset [1,N]$, called the {\em{Cauchon diagram}} of $I$, 
and a sequence of nonzero elements 
\begin{equation}
\label{Yk}
Y_k \in \Fract (R/I), \quad k \in [1,N] \backslash \DD(I),
\end{equation}
which, together with their inverses, generate a copy of a quantum torus inside $\Fract (R/I)$
with commutation relations $Y_k Y_l = \la_{kl} Y_l Y_k$, 
$k,l \in [1,N] \backslash \DD(I)$. This quantum torus contains $R/I$, and is a localization of $R/I$.
\\
We will call the elements \eqref{Yk} the {\em{Cauchon generators}} of $R/I$,
and will denote the set of them by $\CG(R/I)$.   

The sets $\DD(I)$ and $\CG(R/I)$ are defined recursively as follows:

{\em{Case 1,}} $x_N \in I.$ In this case $I = I \cap R_{N-1} + R_N x_N$ and 
we have the isomorphism 
\begin{equation}
\label{phi-case1}
\vp \colon R_N /I \cong R_{N-1}/ (I \cap R_{N-1}),
\end{equation}
the inverse of which is induced by the embedding $R_{N-1} \hra R_N$. 
One sets 
\[
\DD(I) := \DD(I \cap R_{N-1}) \sqcup \{N\} \; \; \mbox{and} \; \; 
\CG(R_N/I):= \vp^{-1} \CG(R_{N-1}/(I \cap R_{N-1}))
\]
and continues recursively with the $\HH$-prime ideal $I \cap R_{N-1}$ 
of $R_{N-1}$. 

{\em{Case 2,}} $x_N \notin I.$ In this case $I[x_N^{-1}] = \oplus_{m \in \Zset} \theta_{x_N}(\lt_{x_N}(I)) x_N^m$, 
\cite[Proposition 2.5(i)]{GeY},  
\begin{equation}
\label{delete-der2}
I':= I[x_N^{-1}] \cap R'_N = \theta_{x_N}(\lt_{x_N}(I))[x_N; \sig_N] \; \; 
\mbox{and} \; \; 
I[x_N^{-1}] = I'[x_N^{-1}].
\end{equation}
We have the isomorphisms 
\begin{multline}
\label{phi-case2}
\vp \colon (R_N/I) [x_N^{-1}] \cong (R'_N[x_N^{-1}])/(I'[x_N^{-1}]) \cong 
\\
\cong
(R'_{N-1}/\theta_{x_N}(\lt_{x_N}(I)))[x_N^{\pm 1}; \sig_N] 
\cong (R_{N-1}/\lt_{x_N}(I)) [ x_N^{\pm 1}; \sig_N]
\end{multline}
where the last map is obtained by applying $\theta_{x_N}^{-1}$ and keeping $x_N$ 
fixed. The first three maps are the canonical isomorphisms induced by \eqref{delete-der1} 
and \eqref{delete-der2}. One sets 
\[
\DD(I) := \DD(\lt_{x_N}(I)) \; \; \mbox{and} \; \; 
\CG(R_N/I):= \vp^{-1} \CG(R_{N-1}/\lt_{x_N}(I)) \sqcup \{ x_N \}
\]
and continues recursively with the $\HH$-prime ideal $\lt_{x_N}(I)$ 
of $R_{N-1}$. 
%%%%%%%%%%%%%%%%%%%%%%%
\subsection{Torus invariant prime ideals of the quantum Schubert cell algebras}
\label{2.4}
The algebra $\UU_q(\g)$ is $Q$-graded by setting
$\deg E_i = \al_i$, $\deg F_i = - \al_i$, and
$\deg K_i^{\pm 1} = 0$. The graded component of 
$\UU_q(\g)$ of degree $\ga \in Q$ will be denoted by $\UU_q(\g)_\ga$.  
The rational character lattice of the $\KK$-torus
\[
\HH:= (\KK^*)^{r}
\]
is identified with the weight lattice $P$ of $\g$ by mapping $\nu \in P$ 
to the character
\[
(t_1, \ldots, t_r) \mt
(t_1, \ldots, t_r)^\nu:= \prod_{i=1}^r t_i^{ \lcor \nu, \al_i\spcheck \rcor }, 
\quad
\forall t_1, \ldots, t_r \in \KK^*.
\] 
The torus $\HH$ acts rationally on $\UU_q(\g)$ by algebra automorphisms by
\[
h \cdot z = h^\ga z \quad \mbox{for} \quad
z \in \UU_q(\g)_\ga, \ga \in Q.  
\]
This action preserves the subalgebras $\UU^\pm[w]$. There is a unique algebra automorphism 
$\om$ of $\UU_q(\g)$ that satisfies 
\begin{equation}
\label{om}
\om(E_i) = F_i, \; \; \om(F_i) = E_i, \; \;  
\om(K_i) = K_i^{-1} \quad
\forall i \in [1,r].
\end{equation}
This automorphism restricts to
an isomorphism $\om \colon \UU^+[w] \cong \UU^-[w]$, $\forall w \in W$.

The Levendorskii--Soibelman straightening law is the following commutation relation 
in $\UU^-[w]$
\begin{multline}
\label{LS}
F_{\be_j} F_{\be_k} - 
q^{ - \lcor \be_k, \be_j \rcor }
F_{\be_k} F_{\be_j}  \\
= \sum_{ {\bf{n}} = (n_{k+1}, \ldots, n_{j-1}) \in \Nset^{\times (j-k-2)} }
p_{\bf{n}} (F_{\be_{j-1}})^{n_{j-1}} \ldots (F_{\be_{k+1}})^{n_{k+1}},
\; \; p_{\bf{n}} \in \KK,
\end{multline}
for all $k < j$. It follows from \eqref{LS} that 
each algebra $\UU^-[w]$ is a symmetric CGL extension with an original presentation 
of the form
\begin{equation}
\label{Uw-p1}
\UU^-[w] = \KK [F_{\be_1}] [F_{\be_2}; \sig_2, \de_2] \ldots [F_{\be_N}; \sig_N, \de_N]
\end{equation}
and reverse presentation 
\begin{equation}
\label{Uw-p2}
\UU^-[w] = \KK [F_{\be_N}] [F_{\be_{N-1}}; \sig^*_{N-1}, \de^*_{N-1}] \ldots [F_{\be_1}; \sig^*_1, \de^*_1].
\end{equation}
The automorphisms $\sig_k$ and $\sig^*_k$ are given by
\[
\sig_k = (h_k \cdot) \quad \mbox{and} \quad \sig_k = (h_k^{-1} \cdot)
\]
where $h_k$ are the unique elements of $\HH$ such that $h^\ga= q^{\lcor \ga, \be_k \rcor}$ for 
for all $\ga \in P$. The skew derivations $\de_k$ and $\de^*_k$ 
are given by 
\[
\de_k(z) = F_{\be_k} z - (h_k \cdot z) F_{\be_k} \quad \mbox{and} \quad 
\de_k^*(z) = F_{\be_k} z - (h_k^{-1} \cdot z) F_{\be_k}.
\]

The $\HH$-primes of 
$\UU^-[w]$ are classified by the following result \cite[Theorem 3.1]{Y-Glasg}.

\bth{Hprim-Uw} For all symmetrizable Kac--Moody algebras $\g$, Weyl group elements $w$,
base fields $\KK$ and non-roots 
of unity $q \in \KK^*$, the poset of $\HH$-primes of $\UU^-[w]$ is isomorphic to 
$W^{\leq w}$ equipped with the Bruhat order. The ideal corresponding 
to $u \in W$, $u \leq w$ is given by
\[
I_w(u) = \{ \lcor c_{\xi, v_{w \la}} \tau \otimes \id, \RR^{w} \rcor \mid \xi \in V(\la)^*, 
\xi \perp \UU^-_q(\g) v_{u \la}, \la \in P^+\}.
\]
\eth

In \cite{Y-Glasg} this result was formulated for simple finite dimensional Lie algebras $\g$. 
However, all proofs in \cite{Y-Glasg} carry over word-by-word to all symmetrizable Kac--Moody 
algebras $\g$. In all proofs one uses the quantized coordinate ring of the corresponding
Kac--Moody group $G$ instead of the quantized coordinate ring of the connected 
simply connected finite dimensional simple Lie group. 
The former is the subalgebra of the dual Hopf algebra $(\UU_q(\g))^*$ 
consisting of the matrix coefficients of all finitely generated integrable highest 
weight $\UU_q(\g)$-modules.    
\subsection{Quantum Richardson varieties}
\label{2.5} 
For $u \in W$, $u \leq w$, denote the canonical projection
\[
\pp \colon \UU^-[w] \to \UU^-[w]/I_w(u). 
\]
The elements $\{ \pp(\De_{u \la, w \la}) \mid \la \in P^+ \}$ are nonzero normal elements of 
$\UU^-[w]/I_w(u)$ and
\begin{equation}
\label{De-comm}
\pp(\De_{u \la, w \la}) z = q^{- \lcor (w+u) \la, \ga \rcor} z \pp(\De_{u \la, w \la}), 
\quad \forall z \in (\UU^-[w]/I_w(u))_\ga, \ga \in Q,
\end{equation}
see \cite[Theorem 3.1 (b), Eq. (3.1)]{Y-Glasg}. It follows from \cite[Theorem 2.6]{Y-plms} that they satisfy 
\begin{equation}
\label{mult-lala}
\De_{u \la_1, w \la_1} \De_{u \la_2, w \la_2} =
q^{\lcor (w-u) \la_1, u \la_2 \rcor} \De_{u (\la_1 + \la_2), w (\la_1 + \la_2)}, 
\quad \forall \la_1, \la_2 \in P^+.
\end{equation}
Given $u \leq w$, one defines the open Richardson variety 
\[
R_{u,w} := B_+ \cdot w B_+ \cap B_- \cdot u B_+ \subset G/B_+.
\]
Its quantized coordinate ring is defined by 
\[
R_q[R_{u,w}]:= \UU^-[w]/I_w(u) [ \pp(\De_{u \la, w \la})^{-1}, \la \in P^+].
\]
This algebra has a canonical rational form over $\Qset[q^{\pm1}]$ whose 
specialization is isomorphic to the coordinate ring of $R_{u,w}$ in the case when 
$\KK= \Cset$, see e.g. \cite[Sect. 4]{Y-cont}. 
In the finite dimensional case, this is proved in \cite[Sect. 4]{Y-cont}; 
the general case is analogous.

Given an integral weight $\la \in P$, let $\la_1, \la_2 \in P^+$ be such that $\la = \la_1 - \la_2$. 
Denote the localized quantum minors
\begin{equation}
\label{De-local}
\De_{u \la, w \la} := 
q^{- \lcor (w-u) \la_1, u \la_2 \rcor} \De_{u \la_1, w \la_1} \De_{u \la_2, w \la_2}^{-1} 
\in \Fract(\UU^-[w]).
\end{equation}
It follows from \eqref{mult-lala} that this definition does not depend 
on the choice of $\la_1, \la_2 \in P^+$.
%%%%%%%%%%%%%%%%%%%%%%%%%%%%%%%%%%%%%%%%%%%%%%%%%%%%%%%
\sectionnew{Contraction of $\HH$-primes in symmetric CGL extensions}
\label{GenContr}
%%%%%%%%%%%%%%%%%%%%
In this section we prove a very general contraction formula 
for the $\HH$--prime ideals of a CGL extension $R$ with 
the subalgebras in the chain \eqref{chain}. This formula is given 
in terms of the Cauchon diagrams with respect to the reverse 
presentation of $R$ which has to do with the chain \eqref{chain-rev}.
\subsection{Statement of main result}
\label{3.1} Let 
\begin{equation}
\label{CGL}
R = \KK[x_1][x_2; \sig_2, \de_2] \ldots [x_N; \sig_N, \de_N]
\end{equation}
be a symmetric CGL extension as in \deref{sCGL} with a rational 
action of the $\KK$-torus $\HH$.
Consider the reverse CGL extension presentation of $R$
\begin{equation}
\label{revCGL}
R = \KK[x_N][x_{N-1}; \sig^*_{N-1}, \de^*_{N-1}] \ldots [x_1; \sig^*_1, \de^*_1]
\end{equation}
and the reverse CGL extension presentation of $R_{N-1}$
\begin{equation}
\label{revCGLN-1}
R_{N-1} = \KK[x_{N-1}][x_{N-2}; \sig^*_{N-2}, \de^*_{N-2}] \ldots [x_1; \sig^*_1, \de^*_1].
\end{equation}
The maps $\sig^*_k$, $\de^*_k$, $k \in [1,N-2]$ in \eqref{revCGLN-1} 
are restrictions of the corresponding maps in \eqref{revCGL}. We use the same symbols for 
simplicity of the notation.

\bde{rever}
For every $\HH$-prime ideal $I$ of a symmetric CGL extension $R$ as in \eqref{CGL},
we will denote by $\DD_{\rev} (I)$ the Cauchon diagram of $I$ with respect to 
the reverse presentation \eqref{revCGL}, where the indices in $\DD_{\rev} (R)$
are recorded in the same way as those of the generators $x$ (without any change 
of the enumeration).

Given a prime ideal $J$ of the symmetric CGL extension 
$R_{N-1}$, denote its Cauchon diagram with respect to the presentation \eqref{revCGLN-1} by 
$\DD_{\rev}(J)$ with the same convention as the one for the ideals of $R$.
\ede

To clarify the convention in \deref{rever}, we give an example. If $R x_1$ is a prime ideal
of $R$, then $\DD_{\rev} (Rx_1) = \{1 \}$ rather than $\{N\}$.

The next theorem provides a very general contraction statement for $\HH$-primes of symmetric CGL 
extensions. It provides a bridge between the two approaches of Goodearl--Letzter (contractions) \cite{GL} 
and Cauchon (deleting derivations) \cite{Ca1} to the $\HH$-prime ideals of symmetric CGL extensions.

\bth{CGLcontr} Let $R$ be a symmetric CGL extension of length $N$. Assume that $I$ and $J$ 
are two $\HH$-prime ideals of $R$ and $R_{N-1}$, respectively,   
such that 
\[
\DD_{\rev}(J) = \DD_{\rev}(I) \cap [1,N-1].
\]
Then 
\[
J = I \cap R_{N-1}.
\]
\eth
\subsection{Proof of \thref{CGLcontr}}
\label{3.2}
We argue by induction on $N$. First, consider the case $1 \in \DD_{\rev}(J)$. 
The assumption $\DD_{\rev}(I) \cap [1,N-1] = \DD_{\rev}(J)$ implies
$1 \in \DD_{\rev}(I)$. Therefore
\[
I = x_1 R + I \cap R_{[2,N]},
\quad J = x_1 R_{N-1} + J \cap R_{[2, N-1]}.
\]
and
\begin{align*}
I \cap R_{N-1} &= \left( x_1 R + I \cap R_{[2,N]} \right) \cap R_{N-1} 
\\
&= x_1 R_{N-1} + I \cap R_{[2,N]} \cap R_{N-1}. 
\end{align*}
The second equality follows from the fact that 
$\{ x_1^{m_1} \ldots x_N^{m_N} \mid m_1, \ldots, m_N \in \Nset \}$ 
is a basis of $R$.

Consider the CGL extension presentations of $R_{[2,N]}$ and $R_{[2,N-1]}$
obtained from \eqref{revCGL} and \eqref{revCGLN-1} 
by removing the last step of the extensions associated to adjoining $x_1$.
The Cauchon diagrams of the $\HH$-prime ideals $I \cap R_{[2,N]}$ and 
$J \cap R_{[2, N-1]}$ of $R_{[2,N]}$ and $R_{[2,N-1]}$ with 
respect to these presentations are 
\[
\DD_{\rev}(I) \backslash \{ 1 \} \quad 
\mbox{and} \quad \DD_{\rev}(J) \backslash \{1\}.
\]
Note that in the second case the generators of $R_{[2,N]}$ are indexed by 
$[2,N]$ in the definition of the diagram. The inductive assumption implies that 
\[
I \cap R_{[2,N]} \cap R_{N-1} = 
(I \cap R_{[2,N]}) \cap R_{[2,N-1]} = J \cap R_{[2,N-1]}.
\]
Thus,
\[
I \cap R_{N-1} = x_1 R_{N-1} + J \cap R_{[2,N-1]} = J.
\]

Next, we consider the case $1 \notin \DD_{\rev}(J)$, which implies $1 \notin \DD_{\rev}(I)$. 
Hence, 
\[
I = \theta_{x_1} \left( \lt_{x_1} (I) \right)[x_1^{\pm 1}; \tau_1 ] \cap R \quad 
\mbox{and} \quad 
J = \theta_{x_1} \left( \lt_{x_1} (J) \right)[x_1^{\pm 1}; \tau_1 ] \cap R_{N-1},
\]
see Section \ref{2.3}.
We have 
\begin{align*}
I \cap R_{N-1} 
&= \theta_{x_1} \left( \lt_{x_1} (I) \right)[x_1^{\pm 1}; \tau_1 ] \cap R_{N-1} 
\\
&= \theta_{x_1} \left( \lt_{x_1} (I) \cap R_{[2,N-1]} \right)[x_1^{\pm 1}; \tau_1 ] \cap R_{N-1}.
\end{align*}
The second equality is proved by recursively applying the property
\[
\theta_{x_1}(a) - a \in 
\oplus_{m=1}^\infty x_1^{-m} R_{[2,k]}, 
\forall a \in R_{[2,k]}, k \in [2,N]. 
\]
This property follows from the definition of the Cauchon map $\theta_{x_1}$ and \prref{sCGL2}.
In this case $\DD_{\rev}(I)$ and $\DD_{\rev}(J)$ coincide with 
the Cauchon diagrams of the $\HH$-primes $\lt_{x_1}(I)$ and $\lt_{x_1} (J)$ 
of $R_{[2,N]}$ and $R_{[2,N-1]}$ obtained from \eqref{revCGL} and \eqref{revCGLN-1} 
by removing the last step of the extensions.   
The inductive assumption implies
\[
\lt_{x_1} (I) \cap R_{[2,N-1]} = \lt_{x_1} (J)
\]
and thus 
\[
I \cap R_{N-1} = 
\theta_{x_1} \left( \lt_{x_1}(J \right))[x_1^{\pm 1}; \tau_1 ] \cap R_{N-1} = J.
\]
This completes the proof of the theorem.
\qed
%%%%%%%%%%%%%%%%
\sectionnew{Contractions of $\HH$-primes of $\UU^-[w]$ 
and sequences of normal elements}
\label{ContrUw}
In this section we prove an explicit formula for the contractions of all $\HH$-prime
ideals of the quantum Schubert cell algebras $\UU^-[w]$ with the intermediate 
subalgebras corresponding to the presentation \eqref{Uw-p1}. For each such 
ideal, the projection of the chain to the prime factor gives a chain 
of subalgebras of the prime factor. We define an explicit sequence of normal 
elements for each such chain. 
\subsection{Contractions}
\label{4.1}
Throughout the section we
fix a Weyl group element $w \in W$ of length $N$ and a reduced expression 
of $w$ as in \eqref{red-expr}. The subexpressions of the latter 
are parametrized by the subsets $D \subseteq [1,N]$. For such a subset $D$, denote 
\[
s_k^D := 
\begin{cases}
s_{i_k}, &\mbox{if} \; \; k \in D,
\\
1, &\mbox{if} \; \; k \notin D
\end{cases} 
\]
and 
\[
w^D_{\leq k} := s_1^D \ldots s_k^D, \quad
w^D_{\geq k} := s_k^D \ldots s_N^D, \quad
w^D_{[j,k]} := s_j^D \ldots, s_k^D, \quad
w^D := w^D_{\leq N} = w^D_{\geq 1}. 
\]
A subexpression is called right positive (respectively left positive) 
if its index set $D$ satisfies $w^D_{\leq k} \leq w^D_{\leq k} s_{i_{k+1}}^D$ $\forall k \in [1,N-1]$,  
(respectively $w^D_{\geq k} \leq s_{i_{k-1}}^D w^D_{\geq k}$, $\forall k \in [2,N]$).
Deodhar and Marsh--Rietsch \cite{De,MR} proved that 
for each $u \in W$ such that $u \leq w$ there exists a unique right positive subexpression 
of \eqref{red-expr}
with total product equal to $u$, i.e., $w^D =u$. Its index set will be denoted 
by $\RP_w(u)$. When one passes from subexpressions of $w$ to those of $w^{-1}$, 
the sets of left and right positive subexpressions are interchanged. 
Hence, for each $u \in W$ such that $u \leq w$ there exists a unique left positive 
subexpression of \eqref{red-expr} with total product equal to $u$. 
Its index set will be denoted by $\LP_w(u)$.

We will use 
the convention of \deref{rever} for diagrams with respect to
reverse presentations. 

\bth{CD} Let $w$ be a Weyl group element with reduced expression 
\eqref{red-expr}. For all Weyl group elements $u \in W$, $u \leq w$, 
the following hold:

(a) the Cauchon diagram $\DD(I_w(u))$ of the $\HH$-prime ideal $I_w(u)$ of 
$\UU^-[w]$ with respect to the presentation \eqref{Uw-p1}
equals the index set $\LP_w(u)$ of the left positive subexpression 
of \eqref{red-expr} with total product $u$, and 

(b) the Cauchon diagram $\DD_{\rev}(I_w(u))$ of $I_w(u)$ of 
$\UU^-[w]$ with respect to the presentation \eqref{Uw-p2}
equals the index set $\RP_w(u)$ of the right positive subexpression 
of \eqref{red-expr} with total product $u$. 
\eth

The first part of the theorem is \cite[Theorem 1.1]{GeY}. 
The proof of the second part is completely analogous.

The intermediate subalgebras for the direct CGL extension presentation 
\eqref{Uw-p1} of $\UU^-[w]$ are given by 
\[
\UU^-[w]_k = \UU^-[w_{\leq k}], \quad k \in [0,N].
\]
For a Weyl group element $u \in W$ such that $u \leq w$, set for brevity
\[
\vec{u}_{\leq k} := w_{\leq k}^{\RP_w(u)}.
\]
The vector notation is suggestive of the definition of right positive subexpression. 
(Right positive subexpressions of reduced expressions of Weyl group elements are picking up indices 
to the far right of the reduced expression.)
A reverse vector notation will be used in relation to left positive subexpressions 
in Section \ref{ContrUwrev}.

Theorems \ref{tCGLcontr} and \ref{tCD}(b) imply at once the following result
describing all contractions of the $\HH$-prime ideals of $\UU^-[w]$ 
with the intermediate subalgebras for the direct presentation of $\UU^-[w]$:

\bth{Uw-contr} Let $w$ be a Weyl group element with reduced expression 
\eqref{red-expr}. For all Weyl group elements $u \in W$ such that $u \leq w$ and 
$k \in [1,N]$, the contractions of the ideal $I_w(u)$ with the 
subalgebras $\UU[w]_k$ are given by 
\[
I_w(u) \cap \UU^-[w_{\leq k}] = I_{w_{\leq k}}(\vec{u}_{\leq k}).
\]
\eth
%%%%%%%%%%%%%%%%%
\subsection{Sequences of normal elements} 
\label{4.2} 
\bde{seq-normal} For an algebra $B$, by {\em{a sequence of nonzero normal elements}} 
\begin{equation}
\label{norm-seq-def}
\De_1, \ldots, \De_l  
\end{equation}
we will mean a sequence of nonzero elements for which there exists a chain of subalgebras
\begin{equation}
\label{chain-def}
B_1 \subset \ldots \subset B_l = B
\end{equation}
such that for all $k$, $\De_k$ is a normal element of $B_k$. 
We will also say that \eqref{norm-seq-def} is a normal 
sequence for the chain \eqref{chain-def} when it is necessary to emphasize the chain of 
subalgebras in the background.
\ede
Note that in general a sequence of normal elements for an algebra $B$ does not consist
of normal elements of $B$. We will construct quantum clusters for an algebra $B$
from sequences of normal elements by removing intermediate terms $\De_k$ that 
are algebraically dependent on the previous terms.  

As in Section \ref{2.5}, for a pair of Weyl group elements $u \leq w$, we will denote by 
\[
\pp \colon \UU^-[w] \to \UU^-[w]/I_w(u)
\]
the canonical projection.
Consider the chain of subalgebras of the prime factor $\UU^-[w]/I_w(u)$
obtained by projecting the intermediate subalgebras for the extension presentation 
\eqref{Uw-p2} 
\begin{equation}
\label{pp-chain}
\pp(\UU^-[w_{\leq 1}]) \subseteq \pp(\UU^-[w_{\leq 2}]) \subseteq \ldots \subseteq 
\pp(\UU^-[w_{\leq N}]) \cong \UU^-[w]/I_w(u).
\end{equation}
\thref{Uw-contr} implies that 
\begin{equation}
\label{factor-k}
\Ker \pp|_{\UU^-[w_{\leq k}]} = I_{w_{\leq k}}(\vec{u}_{\leq k})
\quad \mbox{and} \quad
\pp(\UU^-[w_{\leq k}]) \cong \UU^-[w_{\leq k }] / I_{w_{\leq k}}(\vec{u}_{\leq k}).
\end{equation}

The next theorem constructs a sequence of normal elements for the chain \eqref{pp-chain}.
\bth{chain-norm} 
Let $w$ be a Weyl group element with reduced expression 
\eqref{red-expr} and $u$ be a Weyl group element such that $u \leq w$.
Then for all $k \in [1,N]$ and $\la \in P^+$, 
$\pp(\De_{\vec{u}_{\leq k} \la , w_{\leq k} \la})$ is a nonzero normal element
of the $k$-algebra $\pp(\UU^-[w_{\leq k}])$ in the chain \eqref{pp-chain} 
and more precisely
\begin{equation}
\label{la-norm}
\pp(\De_{\vec{u}_{\leq k} \la , w_{\leq k} \la})
x = q^{-\lcor (w_{\leq k} +\vec{u}_{\leq k}) \la, \ga \rcor} x \pp(\De_{\vec{u}_{\leq k} \la , w_{\leq k} \la}),
\quad \forall \ga \in Q, x \in \pp(\UU^-[w_{\leq k}])_\ga.
\end{equation}

In particular, 
\begin{equation}
\label{sequence}
\pp(\De_{\vec{u}_{\leq 1} \vpi_{i_1} , w_{\leq 1} \vpi_{i_1}}), 
\pp(\De_{\vec{u}_{\leq 2} \vpi_{i_2} , w_{\leq 2} \vpi_{i_2}}), \ldots, 
\pp(\De_{\vec{u}_{\leq N} \vpi_{i_N} , w_{\leq N} \vpi_{i_N}}) \in \UU^-[w]/I_w(u)
\end{equation}
is a sequence with the property that its $k$-th element is a nonzero normal element 
of the $k$-algebra $\pp(\UU^-[w_{\leq k}])$ in the chain \eqref{pp-chain} 
and
\[
\pp(\De_{\vec{u}_{\leq k} \vpi_{i_k} , w_{\leq k} \vpi_{i_k}})
z = q^{- \lcor (w_{\leq k} + \vec{u}_{\leq k}) \vpi_{i_k}, \ga \rcor} 
z \pp(\De_{\vec{u}_{\leq k} \vpi_{i_k} , w_{\leq k} \vpi_{i_k}}),
\]
for all $k \in [1,N], \ga \in Q, z \in \pp(\UU^-[w_{\leq k}])_\ga$.
\eth
\begin{proof} \thref{Uw-contr} (or equivalently \eqref{factor-k}) implies that to prove 
\eqref{la-norm} it is sufficient to establish it for $k =N$, in which case 
it follows from \eqref{De-comm}.
\end{proof}
Taking into account that 
\[
\pp(\De_{\vec{u}_{\leq k} \vpi_{i_k} , w_{\leq k} \vpi_{i_k}})
\in \pp(\UU^-[w_{\leq k}])_{(w_{\leq k} - \vec{u}_{\leq k}) \vpi_{i_k}} 
\]
for $k \in [1,N]$, we obtain the following:
\bco{seq-De} For all pairs of Weyl group elements $u,w \in W$ with $u \leq w$ 
and reduced expressions of $w$, the sequence of elements \eqref{sequence}
of the prime factor $\UU^-[w]/I_w(u)$ is quasi-commuting, more precisely,
\begin{multline}
\label{q-commDe}
\pp(\De_{\vec{u}_{\leq k} \vpi_{i_k} , w_{\leq k} \vpi_{i_k}})
\pp(\De_{\vec{u}_{\leq j} \vpi_{i_j} , w_{\leq j} \vpi_{i_j}}) =
\\
q^{- \lcor (w_{\leq k} + \vec{u}_{\leq k}) \vpi_{i_k}, 
(w_{\leq j} - \vec{u}_{\leq j}) \vpi_{i_j} \rcor} 
\pp(\De_{\vec{u}_{\leq j} \vpi_{i_j} , w_{\leq j} \vpi_{i_j}})
\pp(\De_{\vec{u}_{\leq k} \vpi_{i_k} , w_{\leq k} \vpi_{i_k}})
\end{multline}
for all $1 \leq j < k \leq N$.
\eco
%%%%%%%%%%%%%%%%
\sectionnew{Sequences of normal elements vs. reverse 
Cauchon generators for prime factors of $\UU^-[w]$}
\label{Rel}
In this section we derive an explicit formula expressing the 
Cauchon generators of the $\HH$-prime factors of $\UU^-[w]$ 
with respect to the reverse presentation \eqref{Uw-p2} in terms of the 
sequences of normal elements associated to the direct presentation 
of $\UU^-[w]$ from \thref{chain-norm}. This formula is of monomial nature and 
the exponents have triangular form. As a consequence, the formula yields 
explicit quantum minor generators of the Cauchon quantum tori for 
all $\HH$-prime factors of $\UU^-[w]$ (for the reverse presentation of $\UU^-[w]$).
In particular, this constructs quantum clusters of the $\HH$-prime factors of 
$\UU^-[w]$ in terms of quantum minors. 
\subsection{Statement of the main result}
\label{5.1}
As in the previous section we fix a Weyl group element $w \in W$ and 
a reduced expression \eqref{red-expr} of $w$.
Let $u$ be another Weyl group element such that $u \leq w$.
Recall from Section \ref{4.1} that $\RP_w(u)$ denotes the index set 
of the right positive subexpression of \eqref{red-expr} whose product is
$u$. By \thref{CD}(b), the Cauchon diagram $\DD_{\rev}(I_w(u))$ of the $\HH$-prime ideal 
$I_w(u)$ of $\UU^-[w]$ for the reverse presentation \eqref{Uw-p2} is $\RP_w(u)$. Thus, 
the Cauchon deleting derivation method applied to the reverse 
presentation \eqref{Uw-p2} of $\UU^-[w]$ defines a sequence of nonzero elements 
\[
Y_{k, \rev} \in \Fract ( \UU^-[w]/I_w(u) ), \quad k \in [1,N] \backslash \RP_w(u)
\]
which, together with their inverses, generate a copy of a quantum torus inside \\
$\Fract ( \UU^-[w]/I_w(u) )$. This 
quantum torus contains $\UU^-[w]/I_w(u)$, and 
is a localization of $\UU^-[w]/I_w(u)$.
Following the convention of \deref{rever}, 
the indices of the $Y$-elements match those of the $x$-elements even though  
the $x$-generators are adjoined in the reverse order.

The expression \eqref{red-expr} gives rise to a partial order on $[1,N]$ given by
\begin{equation}
\label{part-ord}
\mbox{
$j \prec k$ 
if $i_j = i_k$ and $j<k$.}
\end{equation}
Set $j \preceq k$ if $j \prec k$ or $j = k$. Extending the 
notation $\vec{u}_{\leq k}$ from the previous section, set 
\[
\vec{u}_{[j,k]} := w^{\RP_w(u)}_{[j,k]}.
\]

Define the following 
integer matrix of size $(N - | \RP_w(u) |)\times N$ whose 
rows are indexed by the set $[1,N] \backslash \RP_w(u)$: 
\[
a_{jk} = 
\begin{cases} 
0, & \mbox{if} \; \; j>k
\\
\lcor \al_{i_j}\spcheck, \vec{u}_{[j+1,k]}(\vpi_{i_k})\rcor = 
\de_{j \preceq k} - \sum_{j < l \preceq k , l \in \RP_w(u)} \lcor \al_{i_j}\spcheck, \vec{u}_{[j+1, l-1]} (\al_{i_l}) \rcor, 
& \mbox{if} \; \; j \leq k,
\end{cases}
\]
where $\de_{j \prec k} :=1$ if $j \prec k$ and $\de_{j \prec k} := 0$ otherwise.
Recall the definition of the quantum minors \eqref{q-min}. The equality 
in the second case follows from Eq. \eqref{second-ajk} in 
\prref{deg}.

\bth{main} Let $\g$ be a symmetrizable Kac--Moody algebra and $w$ be a Weyl group 
element with reduced expression \eqref{red-expr}. Let $u \in W$, $u \leq w$. 
For all base fields $\KK$ and a non-root of unity $q \in \KK^*$,    
in $\Fract ( \UU^-[w]/I_w(u) )$, 
\[
\pp(\De_{\vec{u}_{\leq k} \vpi_{i_k}, w_{\leq k} \vpi_{i_k} }) = 
\prod_{j \in [1,k] \backslash  \RP_w(u) } 
\frac{(q_{i_j}^{-1} - q_{i_j})^{a_{jk}}}{q_{i_j}^{a_{jk}(a_{jk}-1)/2}}
Y_{j, \rev}^{a_{jk}}, \quad \forall k \in [1,N]
\]
where $\pp \colon \UU^-[w] \to \UU^-[w]/I_w(u)$ is the canonical projection.
The product in the right hand side is taken in decreasing order from left to right.
\eth

The special case of the theorem when $\UU^-[w]$ equals the algebra of quantum matrices is due to Cauchon \cite{Ca2}, 
the case $u=1$ (all $w$ and $\g$) was obtained in \cite{GY}.

\bre{triang} (1) The principal minor of the matrix $(a_{jk})$ of size 
$N - \RP_w(u)$ (whose rows and columns are indexed by $[1,N] \backslash \RP_w(u)$) 
is triangular with ones on the diagonal. Thus it is invertible and its inverse 
\[
(b_{jk})_{j, k \in [1,N] \backslash \RP_w(u)} 
\]
is an integral matrix with the same properties. \thref{main} implies that 
\[
Y_{k,\rev} = \zeta_k \prod_{j \in [1,N] \backslash \RP_w(u)} 
\pp(\De_{\vec{u}_{\leq j} \vpi_{i_j}, w_{\leq j} \vpi_{i_j} })^{b_{jk}}, 
\quad \forall k \in [1,N] \backslash \RP_w(u)
\]
where $\zeta_k \in \KK^*$ can be computed explicitly using 
the $q$-commutation relations between the elements $Y_{k,\rev}$. 
The product in the right hand side can be taken in any order, but since 
the terms $q$-commute (\coref{seq-De}), the scalars $\zeta_k$ depend on the choice of order.

(2) \thref{main} and the first part of the remark imply that 
\[
\pp(\De_{\vec{u}_{\leq k} \vpi_{i_k}, w_{\leq k} \vpi_{i_k} }) = 
\theta_k 
\prod_{j \in [1,k] \backslash  \RP_w(u) } 
\pp(\De_{\vec{u}_{\leq j} \vpi_{i_j}, w_{\leq j} \vpi_{i_j} })^{a'_{jk}}, 
\quad \forall k \in \RP_w(u)
\] 
for some integers $a'_{jk}$ and scalars $\theta_k \in \KK^*$ which can be computed explicitly. 
The integers $a'_{jk}$ have the property that $a'_{jk} = 0$ for $j>k$. 

(3) The matrix $(a_{jk})$ has stronger properties than plain triangularity, 
for example,
\[
a_{jk} = \de_{j \preceq k}, \quad \forall j \in [\max \{ l \preceq k \} \cap \RP_w(u), k], 
k \in [1,N].
\]
\ere

\bco{q-torus} For all symmetrizable Kac--Moody algebras $\g$, pairs of 
Weyl group elements $u \leq w$, base fields $\KK$ and non-roots 
of unity $q \in \KK^*$, the nonzero elements 
\[
\pp( \De_{\vec{u}_{\leq k} \vpi_{i_k}, w_{\leq k} \vpi_{i_k}}) \in 
\UU^-[w]/I_w(u), \quad k \in [1,N] \backslash \RP_w(u)
\]
and their inverses 
generate a copy of a quantum torus inside $\Fract(\UU^-[w]/I_w(u))$
with commutation relations \eqref{q-commDe}.
This quantum torus contains $\UU^-[w]/I_w(u)$, and
is a localization of $\UU^-[w]/I_w(u)$
\eco
By the Cauchon procedure (Section \ref{2.3}) the elements
$\{ Y_{k, \rev}^{ \pm 1} \mid k \in [1,N] \backslash \RP_w(u) \}$ 
generate a quantum torus inside $\Fract(\UU^-[w]/I_w(u))$ and this 
quantum torus contains $\UU^-[w]/I_w(u)$. 
The corollary follows from the fact that, by \thref{main}, 
the elements $\{ \pp( \De_{\vec{u}_{\leq k \vpi_{i_k}, w_{\leq k} \vpi_{i_k}}}) \mid 
k \in [1,N] \backslash \RP_w(u)\}$  
generate the same quantum torus.

Recall that a toric frame (with index set $\II \subseteq \Zset$, $|\II|< \infty|$) 
for an algebra $R$ is a map 
\[
M \colon \Zset^\II \to \Fract(R) 
\]
which satisfies the following conditions:
\begin{itemize}
\item For some multiplicatively skewsymmetric 
group bicharacter $\La \colon \Zset^\II \times \Zset^\II \to \KK^*$,
\[
M(f_1) M(f_2) = \La(f_1, f_2) M(f_1 + f_2), \quad \forall f_1, f_2 \in \Zset^\II
\]
(in particular, for the standard basis $\{e_k \mid k \in \II\}$ of $\Zset^\II$, 
$M(e_k)^{\pm 1}$ generate a quantum torus inside $\Fract(R)$ with commutation relations 
$M(e_k) M(e_j) = \La(e_k, e_j)^2 M(e_j) M(e_k)$), 

\item $M (\Nset^\II) \subset R$, and

\item the quantum torus generated by $M(e_k)^{\pm 1}$, $k \in \II$, 
contains $R$. 
\end{itemize} 
A quantum seed for $R$ is a pair consisting of a toric frame and an 
integral $\II \times \II'$ matrix whose principal
part is skewsymmetrizable and which is compatible with the 
cocycle $\La$ in the sense of \cite[Definition 3.1]{BZ} and \cite[\S 2.3]{GY2}. 
(Here $\II' \subseteq \II$ is a set of exchangeable indices.)
We refer to Berenstein--Zelevinsky \cite{BZ} where these notions were introduced.
The case of algebras over $\Qset(q)$ was considered in \cite{BZ} and the general case 
of arbitrary base fields in \cite{GY2}.

We note that \cite{BZ} defines toric frames for division algebras without requiring the third 
condition above. However, if one has a quantum cluster algebra structure on a given algebra $R$, 
then the third condition is a consequence of the quantum Laurent phenomenon. It was shown
in \cite[Sect. 7]{GY2} that the presence of the third property for a family of frames can be used 
in an essential way for the construction of a quantum cluster algebra structure on $R$. This is the reason 
for making it part of the definition. 

\bco{toric-frame} Let $\g$ be a symmetrizable Kac--Moody  algebra, $u \leq w$ a pair of 
Weyl group elements, $\KK$ a base field, and $q \in \KK^*$ a non-root of unity such that 
$\sqrt{q} \in \KK$. The prime factor $\UU^-[w]/I_w(u)$ admits a toric 
frame $M \colon \Zset^{[1,N] \backslash \RP_w(u)} \to \Fract (\UU^-[w]/I_w(u))$
defined by 
\[
M(e_k):=
\pp( \De_{\vec{u}_{\leq k} \vpi_{i_k}, w_{\leq k} \vpi_{i_k}}), \quad \forall 
k \in [1,N] \backslash \RP_w(u)
\]
with respect to the multiplicatively skewsymetric bicharacter given by 
\[
\La(e_k, e_j) :=
\sqrt{q}^{\; - \lcor (w_{\leq k} + \vec{u}_{\leq k}) \vpi_{i_k}, 
(w_{\leq j} - \vec{u}_{\leq j}) \vpi_{i_j} \rcor}, 
\quad
\forall k > j \in \Zset^{[1,N] \backslash \RP_w(u)}.
\]
The toric frame can be extended to a quantum seed of $\UU^-[w]/I_w(u)$
using Leclerc's matrices \cite[Theorem 4.5 and Corollary 4.4]{Le}.
\eco
The compatibility in the last part of the corollary was established by Leclerc in \cite[Sect. 6]{Le}.

The proof of \thref{main} is given in Sections \ref{5.3} and \ref{5.4}; 
Section \ref{5.2} contains some auxiliary degree results that are needed for the proof. 
\prref{main-gen} in Section \ref{5.5} is a generalization 
of the theorem that is needed for the proof of \thref{main2}. For simplicity 
of the exposition we provide full details of the proof of \thref{main} 
and leave the (analogous) proof of \prref{main-gen} to the reader.
%%%%%%%%%%%%%%%
\subsection{Degree considerations}
\label{5.2}
Given a Weyl group element $w$, consider an expression $w=s_{i_1} \ldots s_{i_N}$ which 
is not necessarily reduced. Let 
\begin{equation}
\label{Du}
\SS:= \{ d_1 < \cdots < d_m \} \subseteq [1,N]
\quad \mbox{and} \quad
u = s_{i_{d_1}} \ldots s_{i_{d_m}}.
\end{equation}
We will use the generalization to this setting 
of the notation $\be_k$, $w_{\leq k}$, $w_{[j,k]}$, and the partial order $\prec$ on $[1,N]$. 
For simplicity of the notation set $u_{\leq k}:=  w^\SS_{\leq k}$, 
$u_{[j,k]}:= w^\SS_{[j,k]}$, and $u_{\varnothing}:=1$.
For $j \in [1,k]$ and $\la \in P$, denote
\begin{equation}
\label{a-la}
a_{jk}(\la) :=
\lcor \al_{i_j}\spcheck, u_{[j+1,k]}(\la)\rcor.
\end{equation}
The matrix in Section \ref{5.1} is a special case of this notation: for $\SS = \RP_w(u)$ and $j \leq k$, $a_{jk} = a_{jk}(\vpi_{i_k})$.

\bpr{deg} Let $w=s_{i_1} \ldots s_{i_N} \in W$, and $\SS \subseteq [1,N]$ 
and $u \in W$ be given by \eqref{Du}. 
For all $k \in [1,N]$ and $\la \in P$,
\begin{equation}
\label{deg2}
(w_{\leq k} - u_{\leq k}) \la = 
- \sum_{j \leq k, j \notin \SS} a_{jk}(\la) \be_j.
\end{equation}
Furthermore, for all $i \in [1,r]$ {\em{(}}where $r$ is the rank of $\g${\em{)}} 
and $j \leq k$,
\[
u_{[j+1,k]} (\vpi_i) = 
\om_i - \sum_{j < l \leq k , i_l = i, l \in \SS}  
\lcor \al_{i_j}\spcheck, u_{[j+1,l-1]}(\al_i) \rcor \be_j.
\label{second-ajk}
\]
In particular,
\begin{equation}
\label{ajk}
a_{jk}(\vpi_i) = \de_{i_j, i} - \sum_{j < l \leq k, i_l =i, l \in \SS} 
\lcor \al_{i_j}\spcheck , u_{[j+1, l-1]} (\al_i) \rcor.
\end{equation}
\epr
%The special case of \eqref{deg2} for $\la = - \al_k$ is 
%\[
%u_{\leq k-1} (\al_{i_k}) = 
%\be_k +  \sum_{j< k, j \notin \SS} 
%\lcor \al_{i_j}\spcheck, u_{[j+1,k-1]}(\al_{i_k}) \rcor \be_j
%\end{equation}
%which is an explicit form of \cite[Lemma 4.2]{Y-Glasg}.

\begin{proof} 
For all $j \in [1,k]$, $j \notin \SS$,
\begin{align*}
w_{\leq j} \big( u_{[j+1,k]}(\la) \big) & = 
w_{\leq j-1} s_{i_j} \big( u_{[j+1,k]}(\la) \big)
\\
&= w_{\leq j-1} \big( u_{[j+1,k]}(\la) \big) - 
\lcor \al_{i_j}\spcheck, u_{[j+1,k]}(\la) \rcor \be_j
\\
&= w_{\leq j-1} \big( u_{[j,k]}(\la) \big) - 
a_{jk}(\la) \be_j.
\end{align*}
Adding these identities, proves \eqref{deg2}.

If the set $\{ n \in [j+1,k] \mid i_n =i , n \in \SS \}$ is empty, then the second identity 
in the proposition is trivial. Otherwise, denote by $l$ the maximal element of the set and 
compute 
\[
u_{[j+1,k]} \vpi_i = u_{[j+1,l]} \vpi_i = u_{[j+1,l-1]}(\vpi_i - \al_i).
\]
Then iterate this with $u_{[j+1,k]}$ replaced by $u_{[j+1,l-1]}$.
\end{proof}
\bco{inner-prod} In the setting of \prref{deg}, if $l \notin \SS$, then 
\begin{equation}
\label{cor-ident}
\sum_{j <l, j \notin \SS} \lcor \be_j, \be_l \rcor a_{jk}(\la) - \sum_{l < j \leq k, j \notin \SS} 
\lcor \be_j, \be_l \rcor a_{jk}(\la) = \lcor (w_{\leq k} + u_{\leq k}) \la, \be_l \rcor. 
\end{equation}
\eco
\begin{proof} Note that the integers $a_{jk}(\la)$ only depend on the expression $s_{i_j} \ldots s_{i_k}$ 
and not on the rest of the expression defining $w$. Let $m \in [1,N]$. Applying \eqref{deg2} 
for the expression $s_{i_m} \ldots s_{i_N}$ of $w_{[m,N]}$ gives
\[
(w_{[m,k]} - u_{[m,k]}) \la = - \sum_{m \leq j \leq k, j \notin \SS} a_{jk}(\la) w_{\leq m-1}^{-1}(\be_j)  
\]
and thus
\[
(w_{\leq k} - w_{\leq m-1} u_{[m,k]}) \la = - \sum_{m \leq j \leq k, j \notin \SS} a_{jk}(\la) \be_j.
\]
Using this identity for $m = l$ and $l+1$ and once again \eqref{deg2}, we obtain 
that the left hand side of \eqref{cor-ident} equals 
\begin{align*}
&\sum_{j \leq k, j \notin \SS} \lcor \be_j, \be_l \rcor a_{jk}(\la)
- \sum_{l \leq j \leq k, j \notin \SS} \lcor \be_j, \be_l \rcor a_{jk}(\la)
- \sum_{l < j \leq k, j \notin \SS} \lcor \be_j, \be_l \rcor a_{jk}(\la)
\\
&= \lcor (w_{\leq k} + u_{\leq k}) \la, \be_l \rcor - 
\lcor (w_{\leq l} u_{[l+1, k]} + w_{\leq l-1} u_{[l,k]}) \la, \be_l \rcor.
\end{align*}
The corollary follows from the fact that the last term vanishes,
\begin{align*}
\lcor (w_{\leq l} u_{[l+1, k]} + w_{\leq l-1} u_{[l,k]}) \la, \be_l \rcor
&= \lcor w_{\leq l-1} (s_{i_l} + 1) u_{[l,k]} \la, w_{\leq l-1} \al_{i_l} \rcor =
\\
&= \lcor u_{[l,k]} \la,  (s_{i_l} + 1) \al_{i_l} \rcor = 0,
\end{align*}
where we used the fact that $l \notin \SS$, thus $u_{[l+1, k]} = u_{[l,k]}$.     
\end{proof}
%%%%%%%%%%%%%%%%%%%%%%%%%%%%%%%%%%%%%%%%%%%%%%%%
\subsection{Proof of \thref{main}, part I}
\label{5.3}
Here we prove that 
\begin{equation}
\label{partI}
\pp(\De_{\vec{u}_{\leq k} \vpi_{i_k}, w_{\leq k} \vpi_{i_k} }) = 
\eta_k 
\prod_{j \in [1,k] \backslash  \RP_w(u) } Y_{j, \rev}^{a_{jk}}, \quad \forall k \in [1,N]
\end{equation}
for some $\eta_k \in \KK^*$. In the next subsection we compute the scalars $\eta_k$ explicitly. 

Denote for simplicity
\[
\ol{\De}_k := \pp( \De_{\vec{u}_{\leq k} \vpi_{i_k}, w_{\leq k} \vpi_k}) \in 
\pp( \UU^-[w_{\leq k}]) \subseteq \pp( \UU^-[w]). 
\]
The last part of \thref{chain-norm} implies that
\[
\ol{\De}_k z = q^{-\lcor (w_{\leq k} + \vec{u}_{\leq k}) \vpi_{i_k}, \ga \rcor} 
z \ol{\De}_k, \quad
\forall z \in \Fract(\pp(\UU^-[w_{\leq k}]))_\ga, \ga \in Q.
\]
It follows from the description of the Cauchon procedure in Section \ref{2.3} 
that 
\[
Y_{l, \rev} \in \Fract(\pp(\UU^-[w_{\leq k}])) \subset \Fract(\pp(\UU^-[w]))
\]
for all $l \leq k$, $l \notin \RP_w(u)$, and that $Y_{l, \rev}$ has the same degree as $F_{\be_l}$.
Therefore,
\[
\ol{\De}_k Y_{l, \rev} = 
q^{\lcor (w_{\leq k} + \vec{u}_{\leq k}) \vpi_{i_k}, \be_l \rcor} 
Y_{l, \rev} \ol{\De}_k, \quad \forall l \leq k, l \notin \RP_w(u).
\]
Since $\{Y_{l, \rev}^{\pm 1} \mid l \in [1,k], l \notin \RP_w(u) \}$ generate a quantum torus 
inside $\Fract(\UU^-[w]/I_w(u))$ with commutation relations 
$Y_{j, \rev} Y_{l, \rev} = q^{ - \lcor \be_j, \be_l \rcor} Y_{l, \rev} Y_{j, \rev}$, for $j>l$, we have 
\[
\left( \prod_{j \in [1,k] \backslash \RP_w(u)} Y_{j, \rev}^{a_{j k}} \right) Y_{l, \rev} = 
q^m Y_{l, \rev} \left( \prod_{j \in [1,k] \backslash \RP_w(u)} Y_{j, \rev}^{a_{j k}} \right)
\]  
where 
\[
m = \sum_{j < l, j \notin \RP_w(u)} \lcor \be_j, \be_l \rcor a_{j k} 
- \sum_{l < j \leq k , j \notin \RP_w(u)} \lcor \be_j, \be_l \rcor a_{j k}.
\]
By \coref{inner-prod}, applied to $\la = \vpi_{i_k}$, we have 
$m = \lcor (w_{\leq k} + \vec{u}_{\leq k}) \vpi_{i_k}, \be_l \rcor$. 
Since the division algebra $\Fract(\pp(\UU^-[w_{\leq k}]))$
is generated by the set $\{Y_{l, \rev} \mid l \leq N, l \notin \RP_w(u) \}$, 
\[
\left( \prod_{j \in [1,k] \backslash \RP_w(u)} Y_{j, \rev}^{a_{j k}} \right)
\ol{\De}_k^{\; -1}
\in
Z(\Fract(\pp(\UU^-[w_{\leq k}])).
\]
Applying Eq. \eqref{deg2} in \prref{deg} for $\la = \vpi_{i_k}$, 
we get that the element is in $Z(\Fract(\pp(\UU^-[w_{\leq k}]))^\HH$.

The Goodearl strong $\HH$-rationality result \cite[Theorem II.6.4]{BG} 
for the $\HH$-prime ideals of CGL extensions 
implies that 
\[
Z(\Fract(\pp(\UU^-[w_{\leq k}])))^\HH \cong Z(\Fract(\UU^-[w_{\leq k}]/I_{w_{\leq k}} (\vec{u}_{\leq k})))^\HH=\KK.
\]
This proves that \eqref{partI} is satisfied for some $\eta_k \in \KK^*$.
%%%%%%%%%%%%%%%%%%%%%%%%%%%%%%%%%%%%%%%%%%%%%%%%
\subsection{Proof of \thref{main}, part II}
\label{5.4}
Here we obtain an explicit formula for the scalars $\eta_k$ in \eqref{partI} and complete the 
proof of \thref{main}. The definition of right positive subexpression 
implies that for all $l \in [1,N]$, 

($\st$) {\em{$\RP_w(u)\cap[l,N]$ is the index set of the right positive subexpression of 
$w_{i_l} \ldots w_{i_N}$ 
with total product $\vec{u}_{[l,N]}$.}} 

Denote the scalar in the right hand side of \eqref{partI} for the triple 
$(w_{[l,N]}, \vec{u}_{[l,N]}, k)$ by $\eta_{w_{[l,N]}, \vec{u}_{[l,N]}, k}$.
Here $k \geq l$.
\thref{main}, now follows by iterating the next lemma and using ($\st$) 
at each step.

\ble{ind} For all Weyl group elements $u \leq w$ and $1\leq k \leq N= \ell(w)$, 
\[
\eta_{w,u,k} = 
\begin{cases}
\eta_{w_{[2,N]}, \vec{u}_{[2,N]}, k}, & \mbox{if} \; \; 1 \in \RP_w(u)
\\
\frac{(q_{i_1}^{-1} - q_{i_1})^{a_{1k}}}{q_{i_1}^{a_{1k}(a_{1k}-1)/2}}
\eta_{w_{[2,N]}, \vec{u}_{[2,N]}, k}, & \mbox{if} \; \; 1 \notin \RP_w(u)
\end{cases}
\]
where we set $\eta_{w_{[2,N]}, \vec{u}_{[2,N]}, 1} :=1$.
\ele
\begin{proof} The last extension in the reverse presentation \eqref{Uw-p2} 
of $\UU^-[w]$ is
\begin{equation}
\label{last-st}
\UU^-[w] = T_{s_{i_1}}(\UU^-[w_{[2,N]}])[F_{\be_1}; \sig_1^*, \de_1^*].
\end{equation}
To analyze the nature of the Cauchon procedure applied to this step coupled with 
the effects of the representation theory of $\UU_q(\g)$, we need to consider 
three cases: (1) $1 \in \RP_w(u)$, 
(2) $1 \notin \RP_w(u)$ and $\vec{u}_{\leq k} (\al_{i_1}) \in Q^+$,
(3) $1 \notin \RP_w(u)$ and $\vec{u}_{\leq k} (\al_{i_1}) \in - Q^+$.

Case (1) $1 \in \RP_w(u)$. Then $\vec{u}_{\leq k}^{\, -1} (\al_{i_1}) \in - Q^+$ 
and 
\[
\lcor \xi_{\vec{u}_{\leq k-1} \vpi_{i_k}}, E_{i_1}^m v \rcor =
\lcor S^{-1}(E_{i_1})^m \xi_{\vec{u}_{\leq k-1} \vpi_{i_k}}, v \rcor =
0, \quad \forall m>0, v \in V(\vpi_{i_k}).
\]
It follows from \eqref{q-min2} that
\[
\De_{\vec{u}_{\leq k} \vpi_{i_k}, w_{\leq k} \vpi_k} = 
T_{i_1} (\De_{\vec{u}_{[2,k]} \vpi_{i_k}, w_{[2,k]} \vpi_k}).
\]
\thref{CD} (b) implies that we are 
in the situation of Case 1 in Section \ref{2.3}. So, 
\[
F_{\be_1} \in I_w(u) \quad \mbox{and} \quad 
I_w(u) \cap T_{i_1} (\UU^-[w_{[2,N]}]) = T_{i_1}(I_{w_{[2,N]}}(u_{[2,N]})).
\]
We have the isomorphism
\[
\UU^-[w]/I_w(u) \cong T_{i_1} (\UU^-[w_{[2,N]}]) /  T_{i_1}(I_{w_{[2,N]}}(u_{[2,N]}))
\]
the reverse of which is induced by the embedding 
$T_{i_1} (\UU^-[w_{[2,N]}]) \hra \UU^-[w]$. 
Under this isomorphism,
\[
\De_{\vec{u}_{\leq k} \vpi_{i_k}, w_{\leq k} \vpi_k}  + I_w(u) \mt
T_{i_1} (\De_{\vec{u}_{[2,k]} \vpi_{i_k}, w_{[2,k]} \vpi_k})
+ T_{i_1}(I_{w_{[2,N]}}(u_{[2,N]}))
\]
which proves the lemma in this case.

Case (2) $1 \notin \RP_w(u)$ and $\vec{u}_{\leq k} (\al_{i_1}) \in Q^+$.
Now $\vec{u}_{\leq k} = \vec{u}_{[2,k]}$, 
$a_{1k} = \lcor \al_{i_1} \spcheck, \vec{u}_{[2,k]} \vpi_{i_k} \rcor \geq 0$
and $v_{\vec{u}_{\leq k} \vpi_{i_1}}$ is a highest weight vector for the $\UU_q(\sl_2)$-subalgebra of $\UU_q(\g)$ 
corresponding to the $i_1$-th root with highest weight $a_{1k} \vpi_{i_1}$. 

We will need the following properties of $\UU_q(\sl_2)$-modules and  
the Lusztig braid group action on them: 
\begin{align}
&T_1 v_{n \vpi_1} = \frac{(-q)^n}{[n]_q!} F_1^n v_{n \vpi_1}, \quad
T_1^{-1} v_{n \vpi_1} = \frac{1}{[n]_q!} F_1^n v_{n \vpi_1} \quad \mbox{and} 
\label{br1}
\\
& E_1^n F_1^n v_{n \vpi_1} = ([n]_q!)^2 v_{n \vpi_1}, \quad \forall n>0, 
\label{br2}
\end{align}
see \cite[Eq. 8.6 (3-7) and Lemma 1.7]{Ja}. From them we obtain
\[
E_{i_1}^{a_{1k}} T^{-1}_{i_1} v_{\vec{u}_{\leq k} \vpi_{i_k}} = 
\frac{1}{[a_{1k}]_{q_{i_1}} \! !} E_{i_1}^{a_{1k}} F_{i_1}^{a_{1k}} v_{\vec{u}_{\leq k} \vpi_{i_k}} =
([a_{1k}]_{q_{i_1}} \! !) v_{\vec{u}_{\leq k} \vpi_{i_k}}.
\]
Thus,
\[
\lcor \xi_{\vec{u}_{\leq k} \vpi_{i_k}}, E_{i_1}^{a_{1k}} T_{i_1}^{-1} v \rcor = 
([a_{1k}]_{q_{i_1}} \! !) \lcor \xi_{\vec{u}_{\leq k} \vpi_{i_k}}, v \rcor 
\]
and 
\[
\lcor \xi_{\vec{u}_{\leq k} \vpi_{i_k}}, E_{i_1}^{a_{1k}+ m } T_{i_1}^{-1} v \rcor =0
\]
for all $v \in V(\vpi_{i_k})$, $m >0$. This and the definition of quantum minors \eqref{Rw}--\eqref{q-min2}
imply that in the extension \eqref{last-st}
\[
\De_{\vec{u}_{\leq k} \vpi_{i_k}, w_{\leq k} \vpi_{i_k}} = 
\frac{(q_{i_1}^{-1} - q_{i_1})^{a_{1k}}}{q_{i_1}^{a_{1k}(a_{1k}-1)/2}}
T_{i_1} (\De_{\vec{u}_{[2,k]} \vpi_{i_k}, w_{[2,k]} \vpi_{i_k}}) F_{\be_1}^{a_{1k}} 
+ \mbox{lower order terms in $F_{\be_1}$}.
\]
Applying the result in Section \ref{5.3}, the fact that 
the leading term of the Cauchon map (from Section \ref{2.3}) is the identity
and the equality $Y_{1, \rev} = F_{\be_1}$, 
proves the lemma in this case.

Case (3) $1 \notin \RP_w(u)$ and $\vec{u}_{\leq k} (\al_{i_1}) \in - Q^+$.
Now 
\[
\vec{u}_{\leq k} = s_{i_1} u' \quad 
\mbox{for some} \quad 
u' \in W, \ell(u') = \ell(\vec{u}_{\leq k}) -1.
\]
The above Weyl group elements satisfy 
\[
u' < \vec{u}_{\leq k} = \vec{u}_{[2,k]} \leq w_{[2,k]}
\]
with respect to the Bruhat order.
Moreover, $a_{1k} = \lcor \al_{i_1} \spcheck, \vec{u}_{[2,k]} \vpi_{i_k} \rcor \leq 0$.
For brevity, set
\[
a:= |a_{1k}|.
\]
Then $a_{1k}= -a$ and $v_{\vec{u}_{\leq k} \vpi_{i_1}}$ is a lowest weight vector 
for the $\UU_q(\sl_2)$-subalgebra of $\UU_q(\g)$ corresponding to the $i_1$-th root 
with lowest weight $- a \vpi_{i_1}$. The highest weight vector of this module is 
$v_{u' \vpi_{i_1}}$. This implies that 
\[
S^{-1}(E_{i_1})^m \xi_{u' \vpi_{i_1}} \in (\UU_q^-(\g) \vpi_{\vec{u}_{\leq k} \vpi_{i_1}})^\perp, \quad
\forall m \in [0, a-1].
\]
So, 
\begin{equation}
\label{inideal}
\lcor c_{S^{-1}(E_{i_1})^m, v_{w_{\leq k} \vpi_k}} \tau \otimes \id , \RR^w \rcor 
\in I_{w_{\leq k }}(\vec{u}_{\leq k}), \quad
\forall m \in [0,a-1].
\end{equation}

Given a linear operator $L$ on a vector space $V$, denote its adjoint by 
$L^* \colon V^* \to V^*$, satisfying $\lcor L^* \xi, v \rcor = \lcor \xi, L v \rcor$,
$\forall v \in V, \xi \in V^*$.
For $m \in [0,a]$, denote
\[
D_m = T_{i_1}\Big(
\lcor c_{(T_{i_1}^{-1})^* S^{-1}(E_i)^m \xi_{u' \vpi_{i_k}}, v_{w_{[2,k]} \vpi_{i_k}}} \tau \otimes \id
\rcor \RR^{w_{[2,k]}} \Big) \in T_{i_1}( \UU^-[w_{[2,k]}]).
\]
The properties \eqref{br1}--\eqref{br2} imply
\[
T_{i_1}^{-1} v_{\vec{u}_{\leq k} \vpi_{i_k}} = T_{i_1}^{-1} T_{i_1}^{-1} v_{u' \vpi_{i_1}} = 
(- q_{i_1})^{-a} v_{u' \vpi_{i_k}}
\]
and
\begin{equation}
\label{2ndeq}
E_{i_1}^a v_{\vec{u}_{\leq k} \vpi_{i_k}} = E_{i_1}^a T_{i_1}^{-1} v_{u' \vpi_{i_1}} = 
([a]_{q_{i_1}}\! !) v_{u' \vpi_{i_1}}.
\end{equation}
It follows from the first equality that $(T_{i_1}^{-1})^* \xi_{\vec{u}_{\leq k} \vpi_{i_k}}= \xi_{u' \vpi_k}$. 
Combining this with the definition of the quantum minors \eqref{q-min}, leads to 
\begin{equation}
\label{D0}
D_0 = (- q_{i_1})^{-a} T_{s_1}( \De_{\vec{u}_{[2,k]} \vpi_{i_k}, w_{[2,k]} \vpi_{i_k}}).
\end{equation}
By \eqref{2ndeq}, $(T_{i_1}^{-1})^* S^{-1}(E_i)^a \xi_{u' \vpi_{i_k}} = ([a]_{q_{i_1}} \! !) \xi_{u' \vpi_{i_1}}$.
This and the fact that $S^{-1}(E_{i_1}) \xi_{\vec{u}_{\leq k} \vpi_{i_1}}$ $=0$ 
imply
\begin{equation}
\label{Da}
D_a = ([a]_{q_{i_1}} \! !) T_{i_1} ( \De_{u' \vpi_{i_k}, w_{[2,k]} \vpi_{i_1}}) = 
([a]_{q_{i_1}} \! !) \De_{\vec{u}_{\leq k} \vpi_{i_k}, w_{\leq k} \vpi_{i_k}}.
\end{equation}

The $R$-matrices \eqref{Rw} satisfy $\RR^{w_{\leq k}} = (T_{s_{i_1}} \RR^{w_{[2,k]}}) \RR^{s_1}$. 
Using this and the above mentioned highest weight property of $v_{u' \vpi_{i_1}}$ with 
respect to the $\UU_q(\sl_2)$-subalgebra of $\UU_q(\g)$ associated to the $i_1$-th root, 
after some computations, we obtain 
\[
\lcor c_{S^{-1}(E_{i_1})^{a-m}, v_{w_{\leq k} \vpi_k}} \tau \otimes \id , \RR^w \rcor 
= \sum_{n=0}^{m} \frac{(q_{i_1}^{-1} - q_{i_1})^n}{q_{i_1}^{n(n-1)/2}[n]_{q_{i_1}} \! !}D_{a-m+n}F_{\be_1}^n, 
\quad \forall m \in [0,a].
\]
The left hand side belongs to the ideal $I_w(u)$ for $m \in [1,a]$. 
By induction on $n=0, \ldots, a$, applying the 
$q$-binomial formula
\begin{align*}
&\sum_{n=0}^m \frac{(-1)^{m-n} q^{(m-n)(m-n-1)/2}}{q^{n(n-1)/2} [m-n]_q! [n]_q!}= 
\frac{(-1)^m q^{m^2/2}}{[m]_q !} \sum_{n=0}^m (-1)^n q^{(m-1)n} 
\begin{bmatrix}
a-m
\\
n
\end{bmatrix}_q =
\\
&= \prod_{n=0}^{m-1} ( 1- q^{2(n+1 -m)}),
\end{align*}
we obtain
\[
D_{a-n} = \frac{(-1)^n q_{i_1}^{n(n-1)/2} (q_{i_1}^{-1} - q_{i_1})^n }{[n]_{q_{i_1}} \! !} D_a F_{\be_1}^n
\mod I_w(u).
\]
Combining this with \eqref{D0} and \eqref{Da}, and taking into account that $-a = a_{1k}$, gives
\[
\De_{\vec{u}_{\leq k} \vpi_{i_k}, w_{\leq k} \vpi_{i_k}} F_{\be_1}^{-a_{1k}} = 
\frac{(q_{i_1}^{-1} - q_{i_1})^{a_{1k}}}{q^{a_{1k}(a_{1k}-1)/2}} 
T_{s_1}( \De_{\vec{u}_{[2,k]} \vpi_{i_k}, w_{[2,k]} \vpi_{i_k}})
\mod I_w(u).
\]
The lemma now follows from the fact that the leading term of the 
Cauchon map is the identity, the result in Section \ref{5.3} and the 
fact that $Y_{1, \rev}= F_{\be_1}$.
\end{proof}
%%%%%%%%%%%%%%%
\subsection{A generalization of \thref{main}}
\label{5.5}
For an integral weight $\la \in P$, define the $([1,N] \backslash \RP_w(u))\times [1,N]$ matrix   
\[
a_{jk}(\la) = 
\begin{cases} 
0, & \mbox{if} \; \; j>k
\\
\lcor \al_{i_j}\spcheck, \vec{u}_{[j+1,k]}(\la)\rcor, 
& \mbox{if} \; \; j \leq k,
\end{cases}
\]
which is a specialization of the notation in \eqref{a-la} to the case $\SS = \RP_w(u)$. The following result 
generalizes \thref{main}. Its proof is analogous and is left to the reader. 

\bpr{main-gen} In the setting of \thref{main}, for all $\la \in P$, 
the localized quantum minors \eqref{De-local} satisfy
\[
\pp(\De_{\vec{u}_{\leq k} \la, w_{\leq k} \la }) = 
\prod_{j \in [1,k] \backslash  \RP_w(u) } 
\frac{(q_{i_j}^{-1} - q_{i_j})^{a_{jk}(\la)}}{q_{i_j}^{a_{jk}(\la)(a_{jk}(\la)-1)/2}}
Y_{j, \rev}^{a_{jk}(\la)}, \quad \forall k \in [1,N]
\]
in $\Fract ( \UU^-[w]/I_w(u) )$. 
The product in the right hand side is taken 
in decreasing order from left to right.
\epr
%%%%%%%%%%%%%%%%
\sectionnew{Quantum twist maps for quantum Schubert cell algebras and Richardson 
varieties}
\label{quant-twist}
In this section we define a quantum twist map 
$\Theta_w \colon \UU^-[w^{-1}] \to \UU^-[w]$ that interchanges the 
direct and reverse presentation of the two algebras. It is a quantum version 
of the Fomin--Zelevinsky twist map \cite{FZ}. This map is an algebra 
antiisomorphism. We furthermore prove that 
it restricts to antiisomorphisms 
\[
\Theta_w \colon \UU^-[w^{-1}]/I_{w^{-1}}(u^{-1}) \to \UU^-[w]/I_w(u) \quad
\mbox{and} \quad
\Theta_w \colon R_q[R_{u^{-1},w^{-1}}] \to R_q[R_{u,w}].
\]
\subsection{The quantum twist maps}
\label{6.1a}
For $w \in W$ consider the algebra antiautomorphism
\begin{equation}
\label{twist}
\Theta_w := T_w \tau S \tau \om \colon \UU_q(\g) \to \UU_q(\g)
\end{equation}
where $S$ is the antipode of $\UU_q(\g)$, $\tau$ is the antiautomorphism of $\UU_q(\g)$ 
defined in \eqref{tau}, and $\om$ is the automorphism of $\UU_q(\g)$ defined in \eqref{om}. 
The repetitive use of the map $\tau$ is needed because of the presence of this map in \thref{Hprim-Uw}.  
One checks that $\om \tau = \tau \om$ and $\om S = S \om$, so $\Theta_w$ is also
given by 
\[
\Theta_w = T_w \om \tau S \tau.
\]
\bpr{Thet} For a reduced expression $w = s_{i_1} \ldots s_{i_N}$ 
consider the reduced expression $w^{-1}= s_{i_N} \ldots s_{i_1}$. 
For all $k \in [1,N]$, the antiisomorphism $\Theta_w$ satisifies 
\begin{equation}
\label{Thetaw}
\Theta_w( T_{i_N} \ldots T_{i_{k+1}}(F_{i_k})) = 
\zeta_{w, k} T_{i_1} \ldots T_{i_{k-1}}(F_{i_k})
\end{equation}
for some $\zeta_{w, k} \in \KK^*$. In particular, $\Theta_w$ restricts to an algebra 
antiisomorphism
\[
\Theta_w \colon \UU^-[w^{-1}] \to \UU^-[w].
\]
\epr
\begin{proof} For $\ga = \sum_i n_i \al_i \in Q$, denote $K_\ga := \prod_i K_i^{n_i} \in \UU_q(\g)$.
The antipode satisfies
\begin{equation}
\label{Stau}
S(x) = \zeta_\ga \tau(x) K_{\ga}^{\mp 1}, \quad \forall x \in \UU_q(\b_\pm)_\ga, \ga \in Q
\end{equation}
for some $\zeta_\ga \in \KK^*$ (\cite[Lemma 2.2]{HK}).
This property and the following compatibility 
property of $\tau$ and the braid group action \cite[Eq. 8.18(6)]{Ja}
\[
\tau (T_w x) = T_{w^{-1}}^{-1} ( \tau (x)), \quad
\forall x \in \UU_q(\g), \; w \in W,
\]
imply
\begin{align*}
\tau S \tau \big( T_{i_N} \ldots T_{i_{k+1}} (F_{i_k}) \big) &= \zeta_1
\tau \big( T_{i_N} \ldots T_{i_{k+1}} (F_{i_k}) \big) K_{\pm s_{i_N} \ldots s_{i_{k+1}}(\al_{i_k})} 
\\
&= \zeta_1 T_{i_N}^{-1} \ldots T_{i_{k+1}}^{-1} (F_{i_k}) K_{\pm s_{i_N} \ldots s_{i_{k+1}}(\al_{i_k})} 
\end{align*}
for some $\zeta_1 \in \KK^*$. Hence, 
\begin{align*}
\Theta_w( T_{i_N} \ldots T_{i_{k+1}}(F_{i_k})) &= 
\zeta_1 (T_{i_1} \ldots T_{i_N} \om) \big( T_{i_N}^{-1} \ldots T_{i_{k+1}}^{-1} (F_{i_k}) 
K_{\pm s_{i_N} \ldots s_{i_{k+1}}(\al_{i_k})} \big) 
\\
&= \zeta_2 T_{i_1} \ldots T_{i_N} \big( T_{i_N}^{-1} \ldots T_{i_{k+1}}^{-1} (E_{i_k}) 
K_{\mp s_{i_N} \ldots s_{i_{k+1}}(\al_{i_k})} \big) 
\\
&= \zeta_2 T_{i_1} \ldots T_{i_k} \big( E_{i_k} K^{\mp 1}_{i_k} \big) 
=  \zeta_3 T_{i_1} \ldots T_{i_{k-1}}( F_{i_k}) 
\end{align*}
for some $\zeta_2, \zeta_3 \in \KK^*$ where in the second equality we used the commutation property 
\cite[Eq. 8.18(5)]{Ja} 
\[
\om T_w(x)  = \zeta T_w \om(x), \forall x \in \UU_q(\g)_\ga, \ga \in Q 
\]
for some $\zeta \in \KK^*$ depending on $\ga$.  
\end{proof}
%%%%%%%%%%%%
\subsection{Properties of the quantum twist maps}
\label{6.2a}
For $u \leq w$, denote the canonical projection
\begin{equation}
\label{pi'}
\pp' \colon \UU^-[w^{-1}] \to \UU^-[w^{-1}]/I_{w^{-1}}(u^{-1}).
\end{equation}
By \cite[Theorem 3.3(i)]{FY}, $\{\De_{u \la, w \la}, \la \in P^+\}$ is an Ore subset of $\UU^-[w]$.
We extend $\pp$ and $\pp'$ to projections 
\[ 
\UU^-[w] [\De_{u \la, w \la}, \la \in P^+]^{-1} \to R_q[R_{u,w}], \; 
\UU^-[w^{-1}] [\De_{u^{-1} \la, w^{-1} \la}, \la \in P^+]^{-1} \to R_q[R_{u^{-1},w^{-1}}].
\]
\bth{twist} The following hold for an arbitrary symmetrizable Kac--Moody algebra $\g$, Weyl group element $w$, 
base field $\KK$, and a non-root of unity $q \in \KK^*$:
 
(a) The quantum twist map $\Theta_w$ restricts to an algebra antiisomorphism 
$\Theta_w \colon \UU^-[w^{-1}] \to \UU^-[w]$ which interchanges
the direct and reverse CGL extension presentations \eqref{Uw-p1}--\eqref{Uw-p2} of the two algebras.

(b) For all $u \in W$, $u \leq w$, 
\[
\Theta_w(I_{u^{-1}}(w^{-1})) = I_u(w).
\]

(c) The algebra antiisomorphism
\begin{equation}
\label{fact-Thw}
\Theta_w \colon \UU^-[w^{-1}]/I_{w^{-1}}(u^{-1}) \to \UU^-[w]/I_w(u)
\end{equation}
induces an antiisomorphism 
\begin{equation}
\label{Ruw-Thw}
\Theta_w \colon R_q[R_{u^{-1}, w^{-1}}] \to R_q[R_{u,w}]
\end{equation}
and 
\begin{equation}
\label{De-Thw}
\Theta_w(\pp'(\De_{u^{-1} \la, w^{-1}\la})) = \zeta_\la \pp(\De_{u (u^{-1} \la), w (u^{-1} \la)}), 
\quad \forall \la \in P
\end{equation}
for some $\zeta_\la \in \KK^*$. In the last equality the notation for localized quantum minors \eqref{De-local}
is used for $\la$ and $u^{-1} \la \in P$, respectively. 
\eth
\begin{proof} Part (a) of the theorem follows at once from \prref{Thet}. Part (b) follows 
from \thref{CD} and the first part. 

(c) The antiisomorphism \eqref{fact-Thw} induces an antiisomorphism
\[
\Theta_w \colon \Fract(\UU^-[w^{-1}]/I_{w^{-1}}(u^{-1})) \to \Fract(\UU^-[w]/I_w(u)).
\] 
We will prove that \eqref{De-Thw} holds in $\Fract(\UU^-[w]/I_w(u))$. This 
implies \eqref{Ruw-Thw} and establishes part (c). 

The definition of the quantum twist map $\Theta_w$ gives that
\begin{equation}
\label{grad-Thw}
\Theta_w(z) \in (\UU^-[w]/I_w(u))_{-w\ga}, \quad \forall 
z \in (\UU^-[w^{-1}]/I_{w^{-1}}(u^{-1}))_\ga,  \ga \in Q.
\end{equation}
Eq. \eqref{De-comm} holds for all $\la \in P$. This equation, the fact that $\Theta_w$ 
is an antiisomorphism and the identity
\[
w(w^{-1} \pm u^{-1}) \la = (u \pm w) u^{-1} \la
\]
imply that $\Theta_w(\pp(\De_{u^{-1} \la, w^{-1}\la})) \pp(\De_{u (u^{-1} \la), w (u^{-1} \la)}))^{-1}$
is in the center of the division ring of 
factions of $\UU^-[w]/I_w(u)$.
Furthermore, this identity and \eqref{grad-Thw} imply 
\[
\Theta_w(\pp'(\De_{u^{-1} \la, w^{-1}\la})) \pp(\De_{u (u^{-1} \la), w (u^{-1} \la)})^{-1} \in 
Z(\Fract(\UU^-[w]/I_w(u)))^\HH.
\]
The Goodearl strong rationality result \cite[Theorem II.6.4]{BG} for the torus invariant 
prime ideals of CGL extensions gives $Z(\Fract(\UU^-[w]/I_w(u)))^\HH = \KK$ 
which completes the proof of \eqref{De-Thw} and the theorem.
\end{proof}
The quantum twist map $\Theta_w$ will be used in an essential way in the proof of 
the Berenstein--Zelevinsky conjecture \cite{BZ} in \cite{GY3}.

The theorem has the following corollary for elements of the form $\pp(\De_{u \la, w \la}) \in R_q[R_{u,w}]$ 
for $\la \in P$ that belong to the subalgebra $\UU^-[w]/I_w(u)$.
\bco{integr} We have
\[
\pp(\De_{u \la, w \la}) \in \UU^-[w]/I_w(u) \quad \mbox{for} \; \; \la \in P^+ \cup u^{-1}(P^+),
\]
\eco 
%%%%%%%%%%%%%%%%%%
\sectionnew{Reverse contractions of $\HH$-primes of $\UU^-[w]$ 
and sequences of reverse normal elements}
\label{ContrUwrev}
In this section we describe the contractions of the $\HH$-prime 
ideals of each of the quantum Schubert cell algebras $\UU^-[w]$ with the intermediate 
subalgebras associated to the reverse presentation \eqref{Uw-p2}.
Using the quantum twist map, we also construct an explicit sequence of normal elements for each of these chains
which we call a sequence of reverse normal elements. 
\subsection{Reverse contractions}
\label{6.1}
As before, $w \in W$ denotes a fixed Weyl group element and we work with a fixed
reduced expression \eqref{red-expr} of it. Denote $w_{<k}: = w _{\leq k-1}$, 
($w_{<1}:=1$) and $w_{\geq k} := w_{[k,N]}$.

The intermediate subalgebras for the reverse presentation 
\eqref{Uw-p2} of $\UU^-[w]$ are given by 
\begin{equation}
\label{Uwkrev}
\UU^-[w]_{k, \rev} = \KK \lcor F_{\be_k}, \ldots, F_{\be_N} \rcor = 
T_{w_{<k}} \UU^-[w_{\geq k}], \quad k \in [1,N].
\end{equation}
For a Weyl group element $u \in W$, $u \leq w$, set
\[
\ola{u}_{[j,k]}:= w_{[j,k]}^{\LP_w(u)}
\quad 
\mbox{and}
\quad
\ola{u}_{\geq k} = \ola{u}_{[k,N]}.
\]
The reverse vector notation is suggestive of the definition of left positive subexpression; the point being that 
left positive subexpressions of reduced expressions are picking up indices 
to the far left of the reduced expression.

The following result describes the contractions of all $\HH$-prime ideals of $\UU^-[w]$ with 
the intermediate subalgebras \eqref{Uwkrev} for the reverse presentation of $\UU^-[w]$. 
It follows from Theorems \ref{tCGLcontr} and \ref{tCD}(a). (One can also use 
\thref{twist} (a)-(b), but this is not really needed at this point.)

\bth{Uw-contr2} For all pairs of Weyl group elements $u \leq w$ and reduced expressions 
\eqref{red-expr} of $w$, the contractions of the ideal $I_w(u)$ with 
the subalgebras $\UU^-[w]_{k, \rev}$ are given by 
\[
I_w(u) \cap T_{w_{<k}} \UU^-[w_{\geq k}] = T_{w_{<k}} \Big( I_{w_{\geq k}}(\ola{u}_{\geq k}) \Big),
\quad \forall k \in [1,N].
\]
\eth
%%%%%%%
\subsection{Sequences of reverse normal elements}
\label{6.2}
Consider the canonical projection
\[
\pp \colon \UU^-[w] \to \UU^-[w]/I_w(u). 
\]
The chain of subalgebras 
\[
\UU^-[w]_{N, \rev} \subset \UU^-[w]_{N-1, \rev} \subset \ldots \subset \UU^-[w]_{1, \rev} = 
\UU^-[w]
\]
gives rise to the chain of subalgebras of the prime quotient
\begin{equation}
\label{chain2}
\pp(T_{w_{<N}} \UU^-[w_{\geq N}]) \subseteq 
\ldots \subseteq
\pp(T_{w_{<2}} \UU^-[w_{\geq 2}]) 
\subseteq \pp(\UU^-[w]) \cong \UU^-[w]/I_w(u).
\end{equation}
By \thref{Uw-contr2}, the $k$-th term in this chain is given by 
\begin{multline}
\label{ktermrev}
\pp(T_{w_{<k}} \UU^-[w_{\geq k}]) \cong 
T_{w_{<k}} \UU^-[w_{\geq k}] / ( T_{w_{<k}} \UU^-[w_{\geq k}] \cap I_w(u)) 
\\
\cong
T_{w_{<k}} \UU^-[w_{\geq k}]/T_{w_{<k}} (I_{w \geq k}(\ola{u}_{\geq k}))
\cong \UU^-[w_{\geq k}]/I_{w \geq k}(\ola{u}_{\geq k}).
\end{multline}
For simplicity of the notation we will write
\[
w_{\geq k}^{\; -1} := (w_{\geq k})^{-1}, \quad
{\ola{u}}^{\; -1}_{\geq k}:= (\ola{u}_{\geq k})^{-1}.
\]
\bth{rev-normseq} Assume the setting of \thref{Uw-contr2}. For all $k \in [1,N]$ and 
$\la \in P^+$, 
$\pp \big( T_{w<k} ( \De_{\ola{u}_{\geq k} ( {\ola{u}}^{\; -1}_{\geq k}) \la, 
w_{\geq k} ({\ola{u}}^{\; -1}_{\geq k}) \la }) \big)$ 
is a nonzero normal element of $\pp(\UU^-[w]_{k, \rev})$, and more precisely  
\begin{multline}
\label{nor0}
\pp\Big(T_{w>k} \De_{\, \ola{u}_{\geq k} (\, {\ola{u}}^{\; -1}_{\geq k} \la) , w_{\geq k} (\, {\ola{u}}^{\; -1}_{\geq k} \la)}\Big)
z 
\\
= q^{- \lcor w_{< k} (w_{\geq k} {\ola{u}}^{\; -1}_{\geq k}+1) \la, \ga \rcor} 
z \pp\Big(T_{w>k} \De_{\, \ola{u}_{\geq k} (\, {\ola{u}}^{\; -1}_{\geq k} \la) , w_{\geq k} (\, {\ola{u}}^{\; -1}_{\geq k} \la)}\Big)
\end{multline}
for all $z \in \pp(\UU^-[w]_{k, \rev})_\ga$, $\ga \in Q$.

The sequence
\begin{equation}
\label{sequence2}
\wt{\De}_k:= \pp\big(T_{w_{<N}} \De_{\, \ola{u}_{\geq k} (\, {\ola{u}}^{\; -1}_{\geq k} \vpi_{i_k}), 
w_{\geq k} (\, {\ola{u}}^{\; -1}_{\geq k} \vpi_{i_k})}\big),  \\
k=N,  \ldots, 1
\end{equation}
has the property that its $k$-th element is a nonzero normal element 
of the $k$-th algebra $\pp(\UU^-[w]_{k, \rev})$ in the chain \eqref{chain2}.
In particular, the elements in the sequence quasi-commute,
\[
\wt{\De}_l \wt{\De}_k =
\\
q^{\lcor (w_{\geq k}^{\; -1} + \, {\ola{u}}^{\; -1}_{\geq k}) \vpi_{i_k}, 
(w_{\geq l}^{\; -1} - \, {\ola{u}}^{\; -1}_{\geq l}) \vpi_{i_l} \rcor} 
\wt{\De}_k \wt{\De}_l
\]
for all $1 \leq k < l \leq N$.
\eth
The theorem follows by applying the quantum twist map to the sequence of normal elements
from \thref{chain-norm} for the algebra $\UU^-[w^{-1}]$, using \thref{twist} and the 
identity
\begin{equation}
\label{Thet-id}
\Theta_{w_1 w_2} = T_{w_1} \Theta_{w_2} \quad
\mbox{for} \; \; w_1, w_2 \in W \; \; \mbox{such that} \; \; \ell(w_1 w_2) =
\ell(w_1) + \ell(w_2).
\end{equation}
The inverses of Weyl group 
elements arise from the application of \thref{twist} (c). There are simpler sequences
of normal elements but they do not have the property proved in the next section 
characterizing the Cauchon generators for $\UU^-[w]/I_w(u)$.

We will call the sequence \eqref{sequence2}, {\em{a sequence of reverse 
normal elements}} for $\UU^-[w]/I_w(u)$. It is a sequence of normal 
elements in the sense of \deref{seq-normal} for the chain of subalgebras \eqref{chain2}.
%%%%%%%%%%%%%%%%
\sectionnew{Sequences of reverse normal elements vs.
Cauchon generators for prime factors of $\UU^-[w]$}
\label{RevNorm-Cauchon}
In this section we use the quantum twist maps to obtain explicit expressions for the Cauchon generators of the $\HH$-prime 
factors of $\UU^-[w]$ with respect to the direct presentation \eqref{Uw-p1} of $\UU^-[w]$
in terms of the sequences of reverse normal elements from the previous section. 
The latter are associated to the reverse presentation \eqref{Uw-p2} of $\UU^-[w]$.
This produces another quantum cluster for each $\HH$-prime 
factor of the algebra $\UU^-[w]$. 

In the next section we show that a recursive combination of the results of this section and Section \ref{Rel}
applied in a recursive fashion to subalgebras of $\UU^-[w]$ can be used to construct whole families of toric frames
for the $\HH$-prime factors of $\UU^-[w]$. 
\subsection{Statement of the main result}
\label{7.1}
As before, $w$ denotes a Weyl group element with a fixed 
reduced expression \eqref{red-expr}, and $u$ is a Weyl group element 
with $u \leq w$. \thref{CD}(a) gives that 
the Cauchon diagram $\DD(I_w(u))$ of the $\HH$-prime ideal 
$I_w(u)$ of $\UU^-[w]$ for the direct presentation \eqref{Uw-p1} equals 
the index set $\LP_w(u)$ of the left positive subexpression of \eqref{red-expr} with product
$u$:
\[
\DD(I_w(u)) = \LP_w(u).
\]
So, the Cauchon deleting derivation method applied to the direct 
presentation \eqref{Uw-p1} of $\UU^-[w]$ defines a sequence of nonzero elements 
\[
Y_k \in \Fract ( \UU^-[w]/I_w(u) ), \quad k \in [1,N] \backslash \LP_w(u).
\]
The elements $\{Y_k^{\pm 1} \mid k \in [1,N] \backslash \LP_w(u) \}$ generate a copy of 
a quantum torus inside
$\Fract ( \UU^-[w]/I_w(u) )$ and this 
quantum torus contains $\UU^-[w]/I_w(u)$.

Recall the partial order $\prec$ on $[1,N]$ from \eqref{part-ord}. 
Consider the following integer matrix of size $(N - | \LP_w(u) |)\times N$ whose 
rows are indexed by the set $[1,N] \backslash \LP_w(u)$: 
\[
b_{lk} = 
\begin{cases} 
0, & \mbox{if} \; \; l<k
\\
1, & \mbox{if} \; \; l = k
\\
\lcor \al_{i_l}\spcheck, {\ola{u}}^{\; -1}_{[k,l-1]} (\vpi_{i_k})\rcor = 
\de_{k \prec l} - \sum_{k \preceq j < l, j \in \LP_w(u)} \lcor \al_{i_l}\spcheck, {\ola{u}}^{\; -1}_{[j+1, l-1]} (\al_{i_j}) \rcor, 
& \mbox{if} \; \; l > k.
\end{cases}
\]
The equality in the third case follows from Eq. \eqref{ajk} in \prref{deg}. As in the previous 
section, we write
\[
{\ola{u}}^{\; -1}_{[k,l]}:= (\, \ola{u}_{[k,l]})^{-1}.
\]

\bth{main2} Let $\g$ be a symmetrizable Kac--Moody algebra and $w$ be a Weyl group 
element with reduced expression \eqref{red-expr}. Let $u \in W$, $u \leq w$. 
For all base fields $\KK$ and a non-root of unity $q \in \KK^*$,    
in $\Fract ( \UU^-[w]/I_w(u) )$, 
\[
\pp \big(T_{w_{<k}} \De_{\, \ola{u}_{\geq k} (\, {\ola{u}}^{\; -1}_{\geq k} \vpi_{i_k}), 
w_{\geq k} (\, {\ola{u}}^{\; -1}_{\geq k} \vpi_{i_k}) }\big) = 
\prod_{l \in [k,N] \backslash  \LP_w(u) } 
\frac{(q_{i_l}^{-1} - q_{i_l})^{b_{lk}}}{q_{i_l}^{b_{lk}(b_{lk}-1)/2}}
Y_l^{b_{lk}}, \quad \forall k \in [1,N]
\]
where $\pp \colon \UU^-[w] \to \UU^-[w]/I_w(u)$ is the canonical projection.
The product in the right hand side is taken in a decreasing order from left to right.
\eth

The case of the theorem when $\UU^-[w]$ equals the algebra of quantum matrices was proved by Cauchon 
in \cite{Ca2}, the case $u=1$ (all $w$ and $\g$) was obtained in \cite{GeY}.

\bre{trinag2} Up to a reordering of rows and columns, 
the matrix in \thref{main2} equals the one in \thref{main} for the Weyl group 
elements $u^{-1}$ and $w^{-1}$ (with the reversed to \eqref{red-expr} reduced expression). 
Because of this the matrix in \thref{main2} has the triangular properties 
in \reref{triang} (after reordering of rows and columns).

The formulas in \thref{main2} prove that the quantum torus inside $\Fract(\UU^-[w]/I_w(u))$
generated by 
\[
\{Y_k^{\pm 1} \mid k \in [1,N] \backslash \LP_w(u) \}
\]
also has generators 
\begin{equation}
\label{De-gen}
\big\{ \pp \big(T_{w_{<k}} \De_{\, \ola{u}_{\geq k} (\, {\ola{u}}^{\; -1}_{\geq k} \vpi_{i_k}), 
w_{\geq k} (\, {\ola{u}}^{\; -1}_{\geq k} \vpi_{i_k}) }\big) \mid
k \in [1,N] \backslash \LP_w(u) \big\}.
\end{equation}
This quantum torus contains $\UU^-[w]/I_w(u)$, and is a localizations of the prime factor.
The elements of the second set are monomials in the elements of the first set 
with exponents given by a triangular integral matrix, and vice versa 
the elements of the first set are monomials in the elements of the second set.
Finally, \thref{main2} also implies that the elements 
\[
\big\{ \pp \big(T_{w_{<k}} \De_{\, \ola{u}_{\geq k} (\, {\ola{u}}^{\; -1}_{\geq k} \vpi_{i_k}), 
w_{\geq k} (\, {\ola{u}}^{\; -1}_{\geq k} \vpi_{i_k}) }\big) \mid
k \in \LP_w(u) \big\}
\]
are mononomials in the elements of the set \eqref{De-gen}. One can easily derive explicit 
formulas for this; we leave the details to the reader.
\ere

\bco{toric-frame2} For every symmetrizable Kac--Moody algebra $\g$, pair of Weyl group elements $u \leq w$, 
base field $\KK$, and a non-root of unity $q \in \KK^*$ such that $\sqrt{q} \in \KK^*$, 
the prime factor $\UU^-[w]/I_w(u)$ has a toric 
frame $M \colon \Zset^{[1,N] \backslash \LP_w(u)} \to \Fract (\UU^-[w]/I_w(u))$
given by 
\[
M(e_k):= \pp \big(T_{w_{<k}} \De_{\, \ola{u}_{\geq k} (\, {\ola{u}}^{\; -1}_{\geq k} \vpi_{i_k}), 
w_{\geq k} (\, {\ola{u}}^{\; -1}_{\geq k} \vpi_{i_k}) }\big).
\]
The corresponding multiplicatively skewsymetric bicharacter is given by 
\[
\La(e_l, e_k) :=
\sqrt{q}^{\lcor (w_{\geq k}^{\; -1} + \, {\ola{u}}^{\; -1}_{\geq k}) \vpi_{i_k}, 
(w_{\geq l}^{\; -1} - \, {\ola{u}}^{\; -1}_{\geq l}) \vpi_{i_l} \rcor}, 
\quad
\forall l > k \in \Zset^{[1,N] \backslash \LP_w(u)}.
\]
The toric frame can be augmented to a quantum seed of $\UU^-[w]/I_w(u)$
using Leclerc's matrices \cite[Theorem 4.5 and Corollary 4.4]{Le}, cf. 
\reref{trinag2}.
\eco
\thref{main2} is proved in the next subsections. For the purposes of 
an induction argument, we establish a stronger result. 
For $\la \in P$ and $l \in [1,N]$, denote 
\[
b_{l}(\la) = \lcor \al_{i_l}\spcheck, \ola{u}_{\geq l}(\la)\rcor. 
\]
Then $b_{lk} = b_l \big( \, {\ola{u}}^{\; -1}_{\geq k} \vpi_{i_k} \big)$. 

\bpr{main-gen2} In the setting of \thref{main2}, 
for all $\la \in P$, we have 
\[
\pp(T_{w_{<k}} \De_{\, \ola{u}_{\geq k} \la, w_{\geq k} \la }) =
\prod_{l \in [k,N] \backslash  \LP_w(u) } 
\frac{(q_{i_l}^{-1} - q_{i_l})^{b_l(\la)}}{q_{i_l}^{b_l(\la)(b_l(\la)-1)/2}}
Y_l^{b_l(\la)}, \quad \forall k \in [1,N]
\]
in $\Fract ( \UU^-[w]/I_w(u) )$, recall the notation \eqref{De-local} for localized quantum minors.
The product in the right hand side is taken in a decreasing order from left to right.
\epr
%%%%%%%%%%%%%%%%%%%%%%
\subsection{Proof of \prref{main-gen2}} 
\label{7.2}
For $\la \in P$, set
\[
\ol{\De}_{\la, k, \rev} := \pp(T_{w_{<k}} \De_{\, \ola{u}_{\geq k} \la, w_{\geq k} \la }).
\]
The identity \eqref{Thet-id} for the quantum twist maps and \thref{twist} imply that 
\[
\ol{\De}_{\la, k, \rev} = \zeta_{k, w} \Theta_w 
\pp'\big(\De_{\, {\ola{u}}^{\; -1}_{\geq k} (\, \ola{u}_{\geq k} \la), w^{\; -1}_{\geq k} (\, \ola{u}_{\geq k} \la)}\big)
\quad
\mbox{for some} \; \; \zeta_{w,k} \in \KK^*,
\]
recall \eqref{pi'}.
It follows from \thref{main} (applied to the Weyl group elements $u^{-1}$ and $w^{-1}$)
and \thref{twist} (a) that 
\begin{equation}
\label{partI-7}
\ol{\De}_{\la, k, \rev} = 
\zeta_k
\prod_{l \in [k,N] \backslash  \LP_w(u) } 
Y_l^{b_l(\la)}, \quad \forall k \in [1,N]
\end{equation}
for some $\zeta_k \in \KK^*$.

We obtain an explicit formula for the scalars in \eqref{partI-7} by induction.
The arguments for the inductive statement are different from those in Section \ref{5.4}.
It follows from the definition of 
left positive subexpressions that for all $j \in [1,N]$, 

($\st\st$) {\em{$\LP_w(u)\cap[1,j]$ is the index set of the left positive subexpression of 
$w_{i_1} \ldots w_{i_j}$ with total product $\vec{u}_{\leq j}$.}} 

For $k \leq j$, let $\zeta_{\la, w_{\leq j}, \vec{u}_{\leq j}, k}$ be the scalar in the right hand side 
of \eqref{partI-7} for the quadruple $(\la, w_{\leq j}, \vec{u}_{\leq j}, k)$ and the above choices of reduced 
expressions of $w_{\leq j}$.
\prref{main-gen2} follows by induction 
from the next lemma and ($\st\st$).

\ble{ind2} For all Weyl group elements $u \leq w$ and $1\leq k \leq N= \ell(w)$, 
\[
\zeta_{\la,w,u,k} = 
\begin{cases}
\zeta_{s_{i_N} \la, w_{\leq N-1}, \ola{u}_{\leq N-1}, k}, & \mbox{if} \; \; N \in \LP_w(u)
\\
\frac{(q_{i_N}^{-1} - q_{i_N})^{b_N(\la)}}{q_{i_N}^{b_N(\la)(b_N(\la)-1)/2}}
\zeta_{\la, w_{\leq N-1}, \ola{u}_{\leq N-1}, k}, & \mbox{if} \; \; N \notin \LP_w(u)
\end{cases}
\]
where $\zeta_{\la, w_{\leq N-1}, \vec{u}_{\leq N-1}, N} :=1$.
\ele

Before we proceed with the proof of \leref{ind2}, we establish a 
general fact on Weyl group invariance of localized quantum minors \eqref{De-local} which is 
of independent interest. This fact will also play a role in the next section in connection 
to toric frames for the quantum Richardson varieties.

\bpr{proj} Let $\la \in P$, $w$ be a Weyl group element with reduced expression \eqref{red-expr} 
and $u \leq w$ be such that $N \in \RP_w(u)$. Then 
\begin{equation}
\label{De-le}
\pp( \De_{u (\la), w(\la)}) = 
\pp( \De_{\vec{u}_{\leq N-1} (s_{i_N} \la), w_{\leq N-1} (s_{i_N} \la)})
\end{equation}
where the localized quantum minors use the notation from \eqref{De-local} with 
$\la$ and $s_i(\la)\in P$, respectively.
\epr
\begin{proof} By \thref{Uw-contr}
\[
I_w(u) \cap \UU^-[w_{\leq N-1}] = I_{w_{\leq N-1}} (\vec{u}_{\leq N-1}).
\]
The embedding 
\[
\UU^-[w_{\leq N-1}] \hra \UU^-[w]
\] 
induces the embedding
\[
\vp \colon \UU^-[w_{\leq N-1}]/ I_{w_{\leq N-1}} (\vec{u}_{\leq N-1}) \hra \UU^-[w]/I_w(u).
\]
We will denote by the same letter the extension of this embedding
to the corresponding division rings of fractions. It is easy to see that the Cauchon 
generators of the prime factor on the left with respect to the reverse presentation 
\[
\UU^-[w_{\leq N-1}] = \KK [F_{\be_{N-1}}] [F_{\be_{N-2}}; \sig^*_{N-1}, \de^*_{N-1}] \ldots [F_{\be_1}; \sig^*_1, \de^*_1]
\]
are precisely 
\[
\{\vp^{-1}(Y_{k, \rev}) \mid k \in [1,N-1] \backslash \RP_w(u) \}. 
\]
The equality \eqref{De-le} now follows from the fact that the two sides have the same expressions 
in the sets $\{Y_{k, \rev} \mid k \in [1,N-1] \backslash \RP_w(u) \}$
and 
$\{\vp^{-1}(Y_{k, \rev}) \mid k \in [1,N-1] \backslash \RP_w(u) \}$
given by \prref{main-gen2}. (It is straightforward to see that the exponents 
in the two expressions are the same.)
\end{proof}
\noindent
{\em{Proof of \leref{ind2}}.}
We consider two cases: (1) $N \in \LP_w(u)$ and (2) $N \notin \LP_w(u)$.

Case (1) $N \in \LP_w(u)$. This implies $u s_{i_N} < u$. 
Hence, $N \in \RP_w(u)$ and $\vec{u}_{\leq N-1} = \ola{u}_{\leq N-1}$. 
The lemma now follows from \prref{proj}.

Case (2) $N \notin \LP_w(u)$. We prove the statement of the lemma for $\la \in P^+$. 
The general case follows from the commutation relations between the elements $Y_l$
and the definition of the localized quantum minors \eqref{De-local}. 

The end of the direct presentation \eqref{Uw-p1} 
of $\UU^-[w]$ is 
\begin{equation}
\label{last-stepp}
\UU^-[w] = \UU^-[w_{\leq N-1}][F_{\be_N}; \sig_N, \de_N].
\end{equation}
In this case we are in the situation of case 2 in Section \ref{2.3} and $\ola{u}_{[k, N-1]} = \ola{u}_{\geq k}$.
Using the fact that $v_{w_{\geq k} \la}$ is a lowest weight vector for 
the $\UU_q(\sl_2)$-subalgebra of $\UU_q(\g)$ spanned by $T_{w_{[k, N-1]}} \{ E_{i_N}, F_{i_N}, K_{i_N}^{\pm 1} \}$
with lowest weight $- b_N(\la) \vpi_{i_N} = - \lcor \la, \al_{i_N}\spcheck \rcor \vpi_{i_N}$,
one easily obtains that, with respect to the presentation \eqref{last-stepp}, the leading term of 
\[
T_{w_{<k}} \De_{\, \ola{u}_{\geq k} \la, w_{\geq k} \la } \quad \mbox{is} \quad
\frac{(q_{i_N}^{-1} - q_{i_N})^{b_N(\la)}}{q_{i_N}^{b_N(\la)(b_N(\la)-1)/2}} F_{\be_N}^{b_N(\la)} 
\big( T_{w_{<k}} \De_{\, \ola{u}_{[k, N-1]} \la, w_{[k, N-1]} \la } \big).
\] 
Now the lemma follows from this, \eqref{partI-7} and the fact that the leading term of the Cauchon map 
from Section \ref{2.3} is the identity.
\qed
%%%%%%%%%%
\subsection{A second Weyl group invariance of localized quantum minors}
\label{7a.3}
Analogously to the proof of \prref{proj} one derives the following mirror version of 
it using \prref{main-gen2}. This fact will be needed in the next section for the construction of 
toric frames for the quantum Richardson varieties.

\bpr{proj2} Let $\la \in P$, $w$ be a Weyl group element with reduced expression \eqref{red-expr} 
and $u \leq w$ be such that $1 \in \RP_w(u)$. Then 
\[
\pp( \De_{u\la, w\la}) = 
\pp( T_{s_1} \De_{\vec{u}_{\geq 1} \la, w_{\geq 1}\la})
\]
in the notation from \eqref{De-local}.
\epr
%%%%%%%%%%%%%%%%%%%%%%%%%%%%%
\sectionnew{Families of toric frames for $\Fract(\UU^-[w]/I_w(u))$}
\label{Xi-toricrames}
In this section we construct families of toric frames
for the $\HH$-prime factors of $\UU^-[w]$ and the quantum Richardson algebras 
$R_q[R_{u,w}]$. This is done by a recursive application of the 
results of Sections \ref{Rel} and \ref{RevNorm-Cauchon} to different chains of subalgebras 
of $\UU^-[w]/I_w(u)$. 
%%%%%%%%%%%%%%%%%%
\subsection{Families of chains of subalgebras of $\UU^-[w]$ and contractions of prime ideals}
\label{8.1}
Let $\Xi_N$ be the subset of the symmetric group $S_N$ which consists of all 
permutations $\pi \in S_N$ such that 
\[
\pi(k) = \max \, \pi( [1,k-1]) +1 \; \;
\mbox{or} 
\; \; 
\pi(k) = \min \, \pi( [1,k-1]) - 1, 
\; \; \forall k \in [2,N].
\]
The subset $\Xi_N$ can be equivalently described as the set of all 
$\pi \in S_N$ such that $\pi([1,k])$ is an interval for all $k \in [2,N]$. 

Consider a symmetric CGL extension $R$, recall \deref{sCGL}.
Each $\pi \in \Xi_N$ gives rise to a CGL extension presentation \cite[Remark 6.5]{GY} of $R$,
\begin{equation}
\label{pi-present}
R = \KK [x_{\pi(1)}] [x_{\pi(2)}; \sig'_{\pi(2)}, \de'_{\pi(2)}] 
\cdots [x_{\pi(N)}; \sig'_{\pi(N)}, \de'_{\pi(N)}],
\end{equation}
where 
\[
\sig'_{\pi(k)} := \sig_{\pi(k)}, \; \; h'_{\pi(k)}:= h_{\pi(k)} \; \; \mbox{and} \; \;   
\de'_{\pi(k)} := \de_{\pi(k)}, 
\quad \mbox{if} \; \;  
\pi(k) = \max \, \pi( [1,k-1]) +1
\]
and
\[
\sig'_{\pi(k)} := \sig^*_{\pi(k)}, \; \; 
h'_{\pi(k)}:= h_{\pi(k)}^* \; \; \mbox{and} \; \;   
\de'_{\pi(k)} := \de^*_{\pi(k)}, 
\quad \mbox{if} \; \;  
\pi(k) = \min \, \pi( [1,k-1]) -1.
\]
The direct presentation \eqref{CGL} of $R$ 
corresponds to the identity element $\pi=1$ and the reverse presentation \eqref{revCGL}
to $\pi$ being equal to the longest element of $S_N$. It was proved in 
\cite[Theorem 8.2]{GY2} that, under very mild assumptions, each $\pi \in \Xi_N$ gives rise 
to a quantum seed of $R$ and that those seeds are related by mutations. This was used to 
develop a general theory of quantum cluster algebra structures on symmetric CGL extensions in 
\cite{GY2}. In what follows we use the results of Sections \ref{Rel} and \ref{RevNorm-Cauchon}
to construct families of toric frames for all $\HH$-prime factors of the algebras 
$\UU^-[w]$ indexed by the elements of $\Xi_N$. The prime quotients of an algebra usually behave 
in much more complicated fashion than the algebra itself and cluster structures for 
such are more difficult to construct.

Fix a Weyl group element $w$ and a reduced expression \eqref{red-expr} of it.  
Let $\pi \in \Xi_N$ where $N:= \ell(w)$. For $k \in [1,N]$, 
define $c(k) \leq d(k) \in [1,N]$ by 
\[
[c(k), d(k)] := \pi([1,k]).
\]
By the definition of $\Xi_N$, 
\begin{equation}
\label{pi-cd}
\pi(k) = c(k) \; \; \mbox{or} \; \; d(k).
\end{equation}
The $k$-th intermediate subalgebra of $R$ with respect to the presentation \eqref{pi-present} is $R_{[c(k), d(k)]}$. 
For the quantum Schubert cell algebra $\UU^-[w]$, 
\begin{equation}
\label{Uint}
\UU^-[w]_{[c(k), d(k)]} = T_{w_{< c(k)}}\UU^-[w_{[c(k),d(k)]}]. 
\end{equation}
We have the direct CGL extension presentation of the algebra $\UU^-[w]_{[c(k), d(k)]}$
\[
\UU^-[w]_{[c(k), d(k)]}
= \KK [F_{\be_{c(k)}}] [F_{\be_{c(k)+1}}; \sig_{c(k)+1}, \de_{c(k)+1}] 
\ldots [F_{\be_{d(k)}}; \sig_{d(k)}, \de_{d(k)}]
\]
and the reverse CGL extension presentation of it
\[
\UU^-[w]_{[c(k), d(k)]} =
\KK [F_{\be_{d(k)}}] [F_{\be_{d(k)-1}}; \sig_{d(k)-1}^*, \de_{d(k)-1}^*] 
\ldots [F_{\be_{c(k)}}; \sig_{c(k)}^*, \de_{c(k)}^*].
\]
The automorphisms $\sig_j$, $\sig_j^*$ and skew derivations 
$\de_j$, $\de_j^*$ are the ones from \eqref{Uw-p1} and \eqref{Uw-p2}, restricted 
to the appropriate subalgebras. Another way to look at 
these presentations is to take the two CGL extension 
presentations \eqref{Uw-p1}--\eqref{Uw-p2} of $\UU^-[w_{[c(k),d(k)]}]$ associated 
to the reduced expression
\[
w_{[c(k),d(k)]} = s_{i_{c(k)}} \ldots s_{i_{d(k)}}
\]
and to apply the automorphism $T_{w_{< c(k)}}$ to the corresponding Lusztig root vectors,
taking into account \eqref{Uint}.
The index sets of the left and right positive subexpressions of this expression
will be computed as subsets of $[c(k),d(k)]$ (not of $[1, d(k)-c(k)+1]$).

Let $u$ be a Weyl group element such that $u \leq w$.
Next, we describe the projections of the chain of subalgebras 
\begin{equation}
\label{pi-chain}
\UU^-[w]_{[c(1),d(1)]} \subset \UU^-[w]_{[c(2),d(2)]} \subset \ldots 
\subset  \UU^-[w]_{[c(N),d(N)]} = \UU^-[w]
\end{equation}
into each prime factor $\UU^-[w]/I_w(u)$. Define recursively a sequence of Weyl group 
elements 
\[
u(N):=u, u(N-1), \ldots, u(1) \in W
\]
as follows. Recall \eqref{pi-cd}.

Case (1) $\pi(k+1) = d(k+1)$. Set
\[
u(k) = \min (u(k+1) s_{i_{d(k+1)}}, u(k+1))
\]
with respect to the Bruhat order.

Case (2) $\pi(k+1) = c(k+1)$. Set
\[
u(k) = \min (s_{i_{c(k+1)}} u(k+1), u(k+1) ).
\]

Note that the sequence depends on $\pi$. This dependence will not be shown 
explicitly for simplicity of the notation. The sequence can be equivalently 
defined by setting
\[
u(k) := 
\begin{cases} 
u(k+1) s_{i_{d(k+1)}} , & \mbox{if} \; \; 
d(k+1) \in \RP_{w_{[c(k+1), d(k+1)]}}(u(k+1))
\\
u(k+1), & \mbox{if} \; \; 
d(k+1) \notin \RP_{w_{[c(k+1), d(k+1)]}}(u(k+1)).
\end{cases}
\] 
in the first case and 
\[
u(k) := 
\begin{cases} 
s_{i_{c(k+1)}} u(k+1), & \mbox{if} \; \; 
c(k+1) \in \LP_{w_{[c(k+1), d(k+1)]}}(u(k+1))
\\
u(k+1), & \mbox{if} \; \; 
c(k+1) \notin \LP_{w_{[c(k+1), d(k+1)]}}(u(k+1)).
\end{cases}
\]
in the second.

Since $\pi(k+1) = c(k+1)$ or $d(k+1)$, and 
$[c(k), d(k)] = [c(k+1), d(k+1)] \backslash \{ \pi(k+1) \}$, 
each of the extensions 
\[
\UU^-[w]_{[c(k), d(k)]} \subset \UU^-[w]_{[c(k+1), d(k+1)]}
\]
falls within the framework of subalgebras of quantum Schubert cell algebras 
treated in Sections \ref{ContrUw} or \ref{ContrUwrev}, i.e., subalgebras 
obtained by removing the first or last of the Lusztig root vectors.
Recursively applying Theorems 
\ref{tUw-contr} and \ref{tUw-contr2} to the chain of subalgebras 
\eqref{pi-chain}, we obtain the following,

\bco{gen-contr} For all pairs of Weyl group elements $u \leq w$, reduced expressions 
\eqref{red-expr} of $w$, and elements $\pi \in \Xi_N$, the contractions of the ideal $I_w(u)$ with 
the subalgebras $\UU^-[w]_{[c(k), d(k)]}$ are given by 
\[
I_w(u) \cap \UU^-[w]_{[c(k), d(k)]} = T_{w_{<c(k)}} \Big( I_{w_{[c(k), d(k)]}}(u(k)) \Big),
\quad \forall k \in [1,N].
\]
\eco
Therefore, in the framework of \coref{gen-contr}, the images of the subalgebras \eqref{pi-chain} 
under the projection 
$\pp \colon \UU^-[w] \to \UU^-[w]/I_w(u)$ 
are given by 
\begin{align*}
\pp(\UU^-[w]_{[c(k), d(k)]}) &\cong 
\UU^-[w]_{[c(k), d(k)]}/
T_{w_{<c(k)}} \Big( I_{w_{[c(k), d(k)]}}(u(k)) \Big)
\\
&\cong 
\UU^-[w_{[c(k), d(k)]}] / I_{w_{[c(k), d(k)]}}(u(k)).
\end{align*}

For an arbitrary $\pi \in \Xi_N$, each of the extensions 
\[
\pp(\UU^-[w]_{[c(k), d(k)]}) \subset \pp( \UU^-[w]_{[c(k+1), d(k+1)]})
\]
falls within the framework of those treated in Theorems \ref{tmain} and \ref{tmain2}. 
We will use those results to construct sequences of normal elements
inside the prime factors $\UU^-[w]/I_w(u)$.
%%%%%%%%%%%%%%%%%%%%%%
\subsection{Families of toric frames for the algebras $\UU^-[w]/I_w(u)$ and $R_q[R_{u,w}]$.}
\label{8.2}
Define the following subset of $[1,N]$,
\begin{align*}
D(\pi):= &\{ d(k) \mid k\in [1,N], \pi(k) = d(k), d(k) \in  \RP_{w_{[c(k), d(k)]}}(u(k)) \}
\\     
&\cup \{ c(k) \mid k\in [1,N], \pi(k) = c(k), c(k) \in \LP_{w_{[c(k), d(k)]}}(u(k)) \}.
\end{align*}
(The dependence of the set $D(\pi)$ on $u$ is not explicitly shown for simplicity 
of the notation.)
The second definition of $u(k)$ implies that 
\[
u(k) = w_{[c(k),d(k)]}^{D(\pi) \cap [c(k),d(k)]}, \quad \forall k \in [1,N]
\]
in the notation of Section \ref{4.1}. For $k \in [1,N]$, define the weights 
\[
\la^\pm_{\pi, k} := 
\begin{cases}
w_{< c(k)}( w_{[c(k), d(k)]} \pm u(k) ) \vpi_{i_{d(k)}}, & \mbox{if} \; \; \pi(k) = d(k)
\\
w_{< c(k)}( w_{[c(k), d(k)]} \pm u(k) ) u(k)^{-1} \vpi_{i_{c(k)}}, & \mbox{if} \; \; \pi(k) = c(k)
\end{cases}
\]
and the sequence of elements
\[
\De_{\pi, k}:= 
\begin{cases}
T_{w_{<c(k)}} \De_{u(k) \vpi_{i_{d(k)}}, w_{[c(k),d(k)]} \vpi_{i_{d(k)}}}, & \mbox{if} \; \; \pi(k) = d(k) \\
T_{w_{<c(k)}} \De_{u(k)(u(k)^{-1}\vpi_{i_{c(k)}}), w_{[c(k),d(k)]} (u(k)^{-1} \vpi_{i_{c(k)}})}, & \mbox{if} \; \; \pi(k) = c(k)
\end{cases}
\]
in the notation \eqref{De-local} for localized minors. 
It follows from \coref{integr} that 
\[
p(\De_{\pi, k}) \in  (\UU^-[w]/I_w(u))_{\la_{\pi,k}^-}, \quad \forall k \in [1,N].
\]
By Theorems \ref{tchain-norm} and \ref{trev-normseq}, $\pp(\De_{\pi,k})$ are nonzero 
normal elements of 
$\pp(\UU^-[w]_{[c(k), d(k)]})$ for all $\la \in P$, and more precisely,  
\[
\pp(\De_{\pi, k}) z = 
q^{ - \lcor \la_{\pi,k}^+, \ga \rcor} 
z \pp(\De_{\pi, k})
\]
for all $z \in \pp(\UU^-[w]_{[c(k), d(k)]})_\ga, \ga \in Q$. In particular, 
\[
\pp(\De_{\pi,1}), \ldots, \pp(\De_{\pi,N})
\]
is a sequence of normal elements for the chain of subalgebras of $\UU^-[w]/I_w(u)$ 
consisting of the images of the intermediate subalgebras \eqref{pi-chain} of $\UU^-[w]$ with respect 
to the presentation \eqref{tau}. Furthermore,
\[
\pp(\De_{\pi,k}) \pp(\De_{\pi,j}) = 
q^{ - \lcor \la_{\pi,k}^+, \la_{\pi,j}^- \rcor }  \pp(\De_{\pi,j}) \pp(\De_{\pi, k}), 
\quad \forall k > j \in [1,N].
\]
\bth{toric-frames} Let $\g$ be a symmetrizable Kac--Moody algebra, $u \leq w$ 
a pair of Weyl group elements, $\KK$ a base field, 
and $q \in \KK^*$ a non-root of unity such that $\sqrt{q} \in \KK^*$. 
For all $\pi \in \Xi_N$, the algebras $\UU^-[w]/I_w(u)$ and $R_q[R_{u,w}]$ have a toric 
frame $M \colon \Zset^{[1,N] \backslash D(\pi)} \to \Fract (\UU^-[w]/I_w(u))$
given by 
\[
M(e_k) = \pp(\De_{\pi, k})
\quad \forall 
k \in [1,N] \backslash D(\pi).
\]
The corresponding multiplicatively skewsymmetric bicharacter is given by 
\[
\La(e_k, e_j) :=
\sqrt{q}^{ \; - \lcor \la_{\pi,k}^+, \la_{\pi,j}^- \rcor }, 
\quad \forall k > j \in [1,N] \backslash D(\pi).
\]
\eth
We have $D(1) = \RP_w(u)$ and $D(w_\ci) = \LP_w(u)$ where $w_\ci$ denotes the longest element
of $\Xi_N$. In the special cases of $\pi = 1$ and $\pi = w_\ci$, 
the toric frames in \thref{toric-frames} recover the ones in Corollaries 
\ref{ctoric-frame} and \ref{ctoric-frame2}.
%%%%%%%%%%
\subsection{Proof of \thref{toric-frames}}
\label{8.3}
Before we proceed with the proof of the theorem we establish two lemmas.

\ble{sect9-1} For all symmetrizable Kac--Moody algebras $\g$, Weyl group elements 
$u \leq w$ and $\pi \in \Sig_N$ (where $N = \ell(w)$), the elements 
$p (\De_{u \la, w \la})$ are Laurent monomials in the quantum torus 
generators $\{ M(e_k) \mid k \in [1,N] \backslash D(\pi) \}$ for all 
integral weights $\la$.
\ele 
The elements in the lemma are precisely the normal elements of $\UU^-[w]/I_w(u)$ for the localization defining 
quantum Richardson varieties in Section \ref{2.5}. The lemma follows from Propositions 
\ref{pproj} and \ref{pproj2} and the identity 
\[
\De_{u \vpi_i, w \vpi_i} = \De_{\vec{u}_{\leq N-1} \vpi_i, w_{\leq N-1} \vpi_i}, 
\quad 
\forall i \in [1,N], i \neq i_N,  
\]
the latter in terms of the reduced expression \eqref{red-expr}. 
\ble{sect9-2} In the setting of the previous lemma and the 
reduced expression \eqref{red-expr}
\begin{align*}
p(\De_{\vec{u}_{\leq N-1} \vpi_{i_N}, w \vpi_{i_N}})
&= (q_{i_N}^{-1} - q_{i_N}) p(F_{\be_N}) p(\De_{\vec{u}_{\leq N-1} \vpi_{i_N}, w_{\leq N-1} \vpi_{i_N}}) - p(x)
\\
&= 
\begin{cases} 
0, & \mbox{if} \; \; N \in \RP_w(u) 
\\
p(\De_{u \vpi_{i_N}, w \vpi_{i_N}}), & \mbox{if} \; \; N \notin \RP_w(u)
\end{cases}
\end{align*}
for some $x \in \UU^-[w_{\leq N-1}]$.
\ele
\begin{proof} The first equality follows from \cite[Proposition 4.7]{GeY}. The second case of the 
second equality is straightforward. In the first case, $N \in \RP_w(u)$ implies 
\[
(\vec{u}_{\leq N-1} - u) \vpi_{i_N} = \vec{u}_{\leq N-1} \al_{i_N} \in Q^+ \backslash \{0\}.
\]
Thus, $\xi_{ \vec{u}_{\leq N-1} \vpi_{i_N}} \perp \UU_q^-(\g) v_{u \vpi_{i_N}}$ and 
$\De_{\vec{u}_{\leq N-1} \vpi_{i_N}, w \vpi_{i_N}} \in I_w(u)$.
\end{proof}
\noindent
{\em{Proof of \thref{toric-frames}.}} The only part of the theorem 
that has not been proved yet is that 
$\UU^-[w]/I_w(u)$ and $R_q[R_{u,w}]$ are subalgebras of the quantum torus generated by 
$\{M(e_k)^{\pm 1} \mid k \in [1,N] \backslash D(\pi) \}$. The second statement follows from the first 
and \leref{sect9-1}. To establish the first statement, we prove by induction 
on $l=1, \ldots, N$, the stronger fact that

{\em{$p(\UU^-[w]_{[c(l),d(l)]})$ is a subalgebra of the quantum subtorus of 
$\Fract (\UU^-[w]/I_w(u))$ generated by $\{M(e_k)^{\pm 1} \mid k \in [1,l] \backslash D(\pi) \}$.}}

Denote by $\TT_l$ the quantum subtorus of 
$\Fract (\UU^-[w]/I_w(u))$ generated by $\{M(e_k)^{\pm 1} \mid k \in [1,l] \backslash D(\pi) \}$.
Assume the validity of the statement for $l-1$. Consider the case $\pi(l) = d(l)$. 
Applying Lemmas \ref{lsect9-1} and \ref{lsect9-2}, we obtain that
\[
p(F_{\be_l}) \in \big( \TT_l - p(\UU^-[w]_{[c(l-1),d(l-1)]}) \big)
p\big(\De_{u(l-1) \vpi_{i_l}, w_{[c(l-1), d(l-1)] \vpi_{i_l}}}\big)^{-1}
\] 
The space on the right is a subspace of $\TT_l$ because of
the inductive assumption and 
the fact that the last term is a Laurent monomial 
in the generators of $\TT_l$ by \leref{sect9-1}. The case $\pi(l) = c(l)$ is handled 
similarly by applying the quantum twist map $\Theta_{w_{[c(l),d(l)]}}$ to the equalities
in \leref{sect9-2}.
\qed
%%%%%%%%%%%%%%%%%%%%%% References %%%%%%%%%%%%%%%%%%%%%%%%%%%%%%%%%%%%%%%

%%%%%%%%%%%%%%%%%%%%%%%%%%%%%%%%%%%%%%%%%%%%%%%%%%%%%%%%%%%%%%%%%%%%%%%%%%%%%%%
%%%%%%%%%%%%%%%%%%%%%%%%%%%%%%%%%%%%%%%%%%%%%%%%%%%%%%%%%%%%%%%%%%%%%%%%%%%%%%
\end{document}